\documentclass{amsart}
\usepackage[square,comma,numbers,sort&compress]{natbib}
\usepackage{fullpage}

\usepackage[linesnumbered,ruled,vlined]{algorithm2e}
\usepackage{tikz}
\usetikzlibrary{shapes}

\usepackage{amssymb}

\usepackage{mathtools}

\usepackage{physics}

\usepackage{todonotes}
\usepackage{cleveref}
\usepackage{longtable}[=v4.13]%

\newcommand{\bu}{\mathbf{u}}
\newcommand{\bv}{\mathbf{v}}
\newcommand{\bk}{\mathbf{k}}

\Crefname{table}{Table}{Tables}
\Crefname{algorithm}{Algorithm}{Algorithms}

\begin{document}




\title[AOT Data Assimilation for MPAS-Ocean]{Application of Continuous Data Assimilation in High-Resolution Ocean Modeling}

\author[Larios A. et.~al.]{Adam Larios, Mark Petersen, and Collin Victor}





\begin{abstract}
We demonstrate a formulation of the Azouani-Olson-Titi (AOT) algorithm in the MPAS-Ocean implementation of the primitive equations of the ocean, presenting global ocean simulations with realistic coastlines and bathymetry.
We observe an exponentially fast decay in the error before reaching a certain error level, which depends on the terms involved and whether the AOT feedback control term was handled implicitly or explicitly. A wide range of errors was observed for both schemes, with the implicit scheme typically exhibiting lower error levels, depending on the specific physical terms included in the model. Several factors seem to be contributing to this wide range, but the vertical mixing term is demonstrated to be an especially problematic term. This study provides insight into the promises and challenges of adapting the AOT algorithm to the setting of high-resolution, realistic ocean models.
\end{abstract}

%
%
%
%
%
\maketitle


\section{Introduction}


\noindent
In the modeling of turbulent flows, such as ocean dynamics, there are many non-trivial issues that arise. One such issue is the initialization of the model. That is, even if one could evolve the model forward in time exactly, one still would require an initial state to evolve forward from. For many dynamical systems of interest it is infeasible to obtain a complete set of initial data with which to initialize the model. One way to mitigate this problem is to use a set of techniques, collectively known as data assimilation, to incorporate observational data into the model.

Data assimilation refers to a class of schemes that combine observational data with a physical model to better predict the future state of a physical system. Classical methods of data assimilation include the Kalman filter and its variants, as well as by variational methods such as 3D/4D Var. For a look at classical methods of data assimilation see, e.g., \cite{Asch_Bocquet_Nodet_2016_DA_book,Law_Stuart_Zygalakis_2015_book} and the references therein. In the past decade, a new data assimilation paradigm known as the Azouani-Olson-Titi (AOT) continuous data assimilation algorithm has developed  \cite{Azouani_Olson_Titi_2014,Azouani_Titi_2014}. 

The AOT algorithm was first developed in the context of the 1D Allen--Cahn Equation in \cite{Azouani_Titi_2014}. It was then given in full generality in  \cite{Azouani_Olson_Titi_2014} for general classes of interpolants in the context of the 2D incompressible Navier--Stokes Equations. 
Since its inception in 2014, the AOT algorithm has been applied successfully to many other dissipative dynamical systems. 
This includes the Cahn--Hilliard Equations \cite{Diegel_Rebholz_2021}, the 3D primitive equations  of the ocean \cite{Pei_2019,Carlson_VanRoekel_Petersen_Godinez_Larios_2021}, the surface quasi-geostrophic equations \cite{Jolly_Martinez_Olson_Titi_2018_blurred_SQG,Jolly_Martinez_Titi_2017}, the Kuramoto--Sivashinsky Equations \cite{Larios_Pei_2017_KSE_DA_NL,Lunasin_Titi_2015,Pachev_Whitehead_McQuarrie_2021concurrent}, and various regularizations of NSE \cite{Albanez_Nussenzveig_Lopes_Titi_2016,Gardner_Larios_Rebholz_Vargun_Zerfas_2020_VVDA,Larios_Pei_2018_NSV_DA} to name a few. This algorithm has also been used in downstaling a general circulation model of the atmosphere \cite{Desamsetti_Dasari_Langodan_Knio_Hoteit_Titi_2019_WRF}.
There have also been numerous studies that utilize the AOT algorithm to recover the solutions to dynamical systems with a variety of restrictions placed on the observational data. This includes observations that include measurement error \cite{Bessaih_Olson_Titi_2015,Foias_Mondaini_Titi_2016,Celik_Olson_2022}, incomplete observations of each variable \cite{Farhat_Jolly_Titi_2015,Farhat_Lunasin_Titi_2016abridged,Farhat_Lunasin_Titi_2016benard,Farhat_Lunasin_Titi_2016_Charney,Farhat_Lunasin_Titi_2017_Horizontal}, observations that are discrete or averaged in time \cite{Larios_Pei_Victor_2023_second_best,Foias_Mondaini_Titi_2016,Celik_Olson_2022,Jolly_Martinez_Olson_Titi_2018_blurred_SQG}. There have also been numerous studies adjusting the gathering and assimilation of observational data, this includes assimilating observational data on a localized patch of the domain \cite{Biswas_Bradshaw_Jolly_2022}, using mobile observers to gather data \cite{Larios_Victor_2019,Franz_Larios_Victor_2022, Biswas_Bradshaw_Jolly_2022}, and using nonlinear variations of the feedback control term 
\cite{Du_Shiue_2021,Carlson_Larios_Titi_2023_nlDA}. 
Related algorithms which use interpolated or filtered data as in AOT, but insert the results directly have been studied in \cite{Celik_Olson_Titi_2019,Hayden_Olson_Titi_2011,Larios_Pei_Victor_2023_second_best,Olson_Titi_2003}.
Recently there has been work done on modifying this algorithm for use not only in recovering solutions, but in system identification. 
In particular, the AOT algorithm has been modified to simultaneously recover the true solution of a dynamical system along with the viscosity \cite{Carlson_Hudson_Larios_2020, Farhat_GlattHoltz_Martinez_McQuarrie_Whitehead_2019} or the forcing \cite{Martinez_2022,Farhat_Larios_Martinez_Whitehead_2023_force}.

The AOT algorithm differs from classical methods of data assimilation by introducing a feedback control term at the PDE level. We note that the AOT algorithm appears similar to a form of data assimilation known as nudging, or Newtonian relaxation proposed and studied in \cite{Anthes_1974_JAS,Hoke_Anthes_1976_MWR}. This similarity is superficial, as the AOT algorithm applies a spatial interpolation which has a large impact on, implementation, practical usage, and convergence rates for the algorithm. For an in depth examination of nudging methods (as opposed to AOT methods) see e.g. \cite{Lakshmivarahan_Lewis_2013} and the references contained within. 

Consider a dynamical system, e.g. the primitive ocean equations, written abstractly as:
\begin{equation}
    \frac{d}{dt}u = F(u).
\end{equation}
Note that, in the above equation, we do not know the initial state $u(0)=u_0$ of the flow. Instead, we observe the velocity at certain spatial locations and incorporate this into the AOT algorithm which is given as follows:
\begin{align}
    \frac{d}{dt} v &= F(v) + \mu \left(I_\delta(u) - I_\delta(v)\right)\\
    v(0,x) &= v_0
\end{align}
Here $\mu>0$ is a constant feedback control term and $I_\delta$ is an interpolation operator corresponding to some characteristic spatial length scale $\delta>0$. Of course, $u$ is unknown, but we suppose that we have access to certain measurements of $u$, for instance, observations of $u$ at finitely many spatial locations.  We then interpolate these observations globally to form $I_\delta(u)$.  The initial data for the reference flow, $u$, is also unknown, but the initial flow for the simulated flow, $v$, is given by $v(0)=v_0$, where $v_0$ is some admissible but otherwise arbitrary initial data. For instance, one may choose $v_0\equiv0$.  
In the present work, for simplicity, we assume the observations are exact and are taken continuously in time.  
However, it has been shown that the AOT algorithm is robust to removing these assumptions.  For example, noisy observations and stochastic forcing have been studied in \cite{Bessaih_Olson_Titi_2015,Carlson_Hudson_Larios_Martinez_Ng_Whitehead_2021}, and sparse-in-time observations have been studied in \cite{Carlson_Hudson_Larios_Martinez_Ng_Whitehead_2021,Foias_Mondaini_Titi_2016,Larios_Pei_Victor_2023_second_best}.  The AOT algorithm also does not need the model parameters to be known perfectly (or at all) in many cases, and indeed, AOT can be used in tandem with other algorithms to recover unknown parameters, as has been studied in \cite{Carlson_Hudson_Larios_Martinez_Ng_Whitehead_2021,Carlson_Hudson_Larios_2020,Larios_Pei_2018_NSV_DA,Pachev_Whitehead_McQuarrie_2021concurrent,Farhat_GlattHoltz_Martinez_McQuarrie_Whitehead_2019}.
We detail below in \Cref{background} our implementation of the AOT algorithm specific to the governing equations of the MPAS-Ocean model.


It was shown in \cite{Pei_2019} that the AOT algorithm, in the context of the 3D primitive equations of the ocean, is globally well-posed and exhibits exponentially fast convergence when applied to the 3D primitive equations with some idealizations.
In this work we examine computationally the performance of the AOT algorithm when applied to a realistic ocean model.
Specifically, we test this algorithm computationally against a high resolution ocean model with realistic coastlines.
We consider in particular the effect of assimilating the tracer quantities, temperature and salinity.
We find that the assimilation of the tracers improves the level of convergence of the kinetic energy of the system, however we do not obtain convergence to the level of machine precision, except in certain cases with simplified physics.
We note that the level of convergence we obtain is comparable to those of recent studies of this algorithm for ocean models including \cite{Carlson_VanRoekel_Petersen_Godinez_Larios_2021} where the simulations were done in the context of an idealized mesoscale eddy test case
and \cite{Hammoud_Titi_Hoteit_Knio_2022}, where the performance of the AOT algorithm was compared to a physics-informed deep neural network in a downscaling model. 

This work is organized as follows. In  \Cref{background} we discuss the formulation of the AOT algorithm equations and the ocean model equations. In \Cref{methods} we describe the experimental setup utilized to test the AOT algorithm, the numerical scheme, details of the MPAS-Ocean model, and we show the results of test utilized to verify the correct implementation and performance of the AOT algorithm.
In \Cref{computational results} we detail the results of our simulations and the effect that varying parameters of the system has on convergence.
We end with some concluding remarks in \Cref{conclusion}.


\begin{figure}
    \centering    
    \includegraphics[width=.4\textwidth, height=4cm]{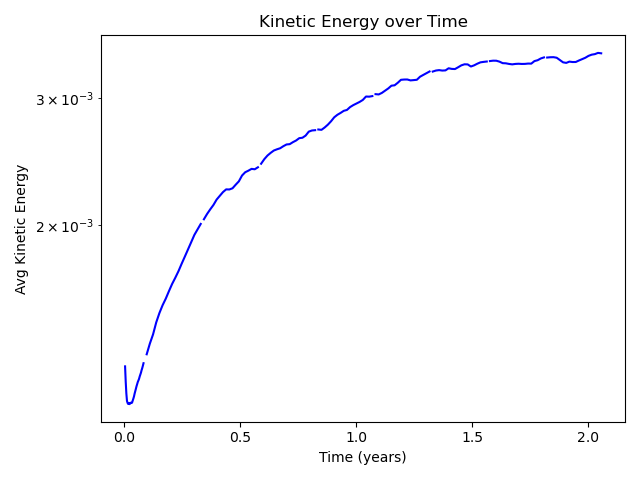}    
    \includegraphics[width=.4\textwidth, height=4cm]{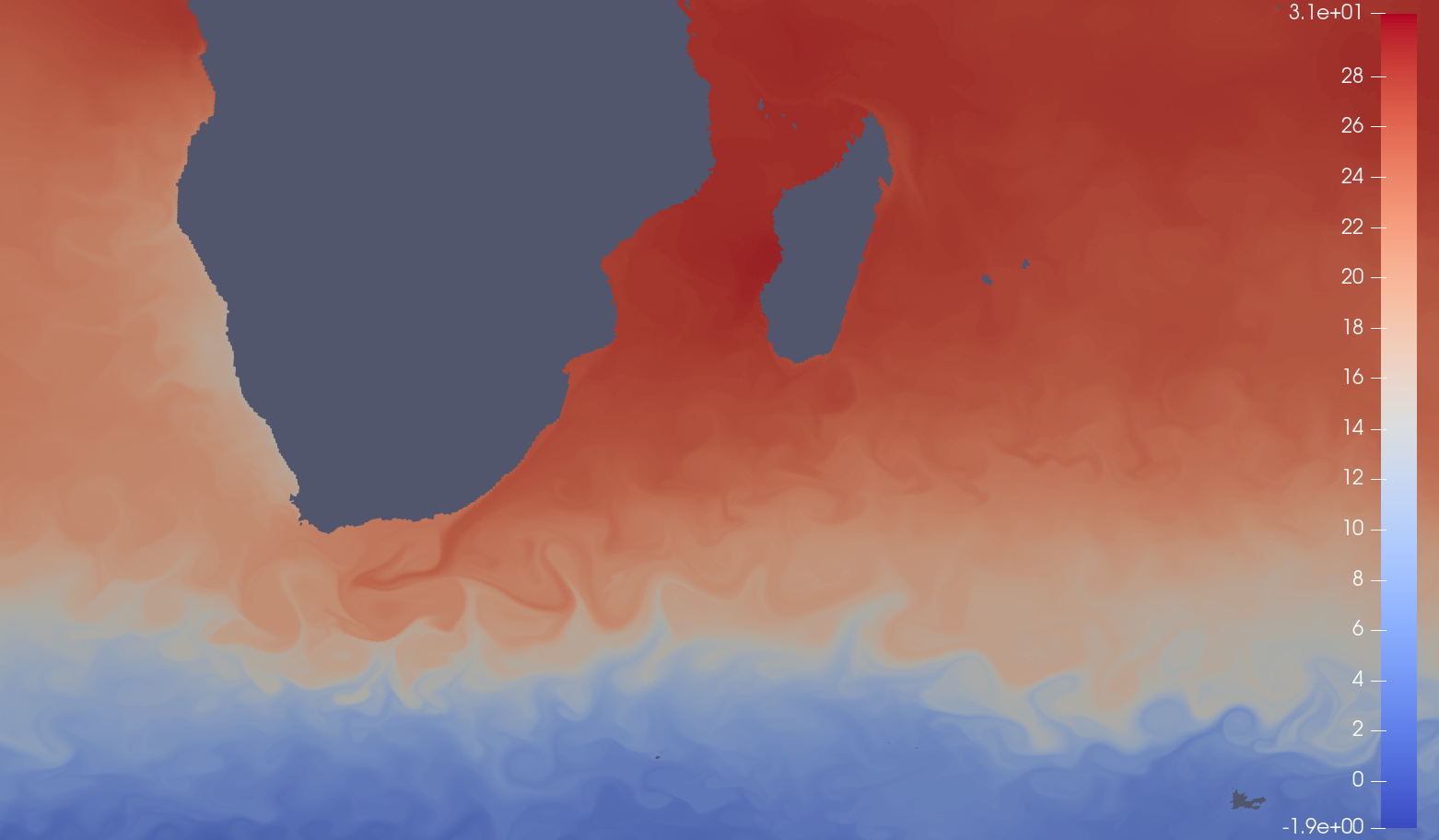}    
    \includegraphics[width=.4\textwidth, height=4cm]{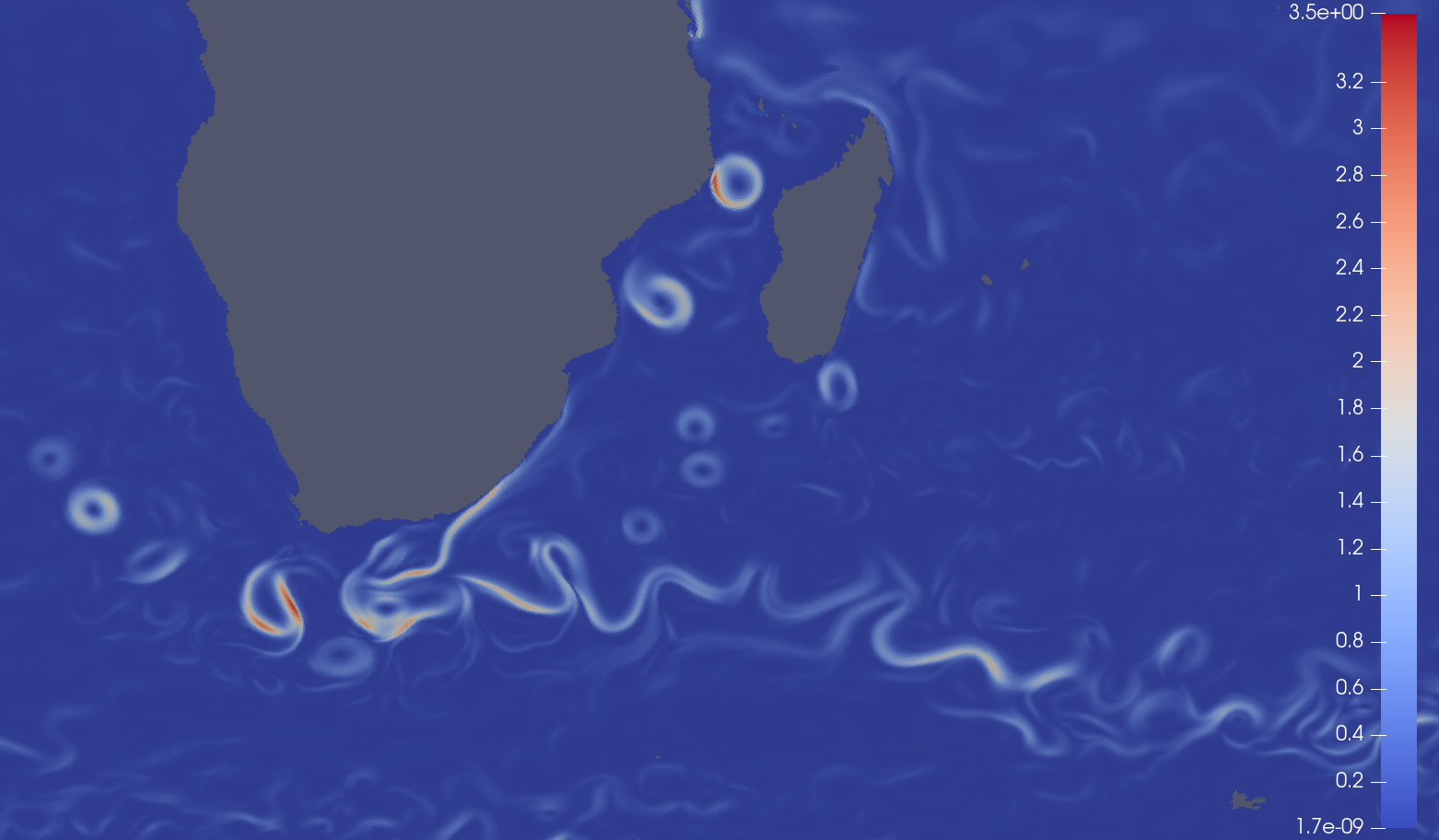}
    \includegraphics[width=.4\textwidth, height=4cm]{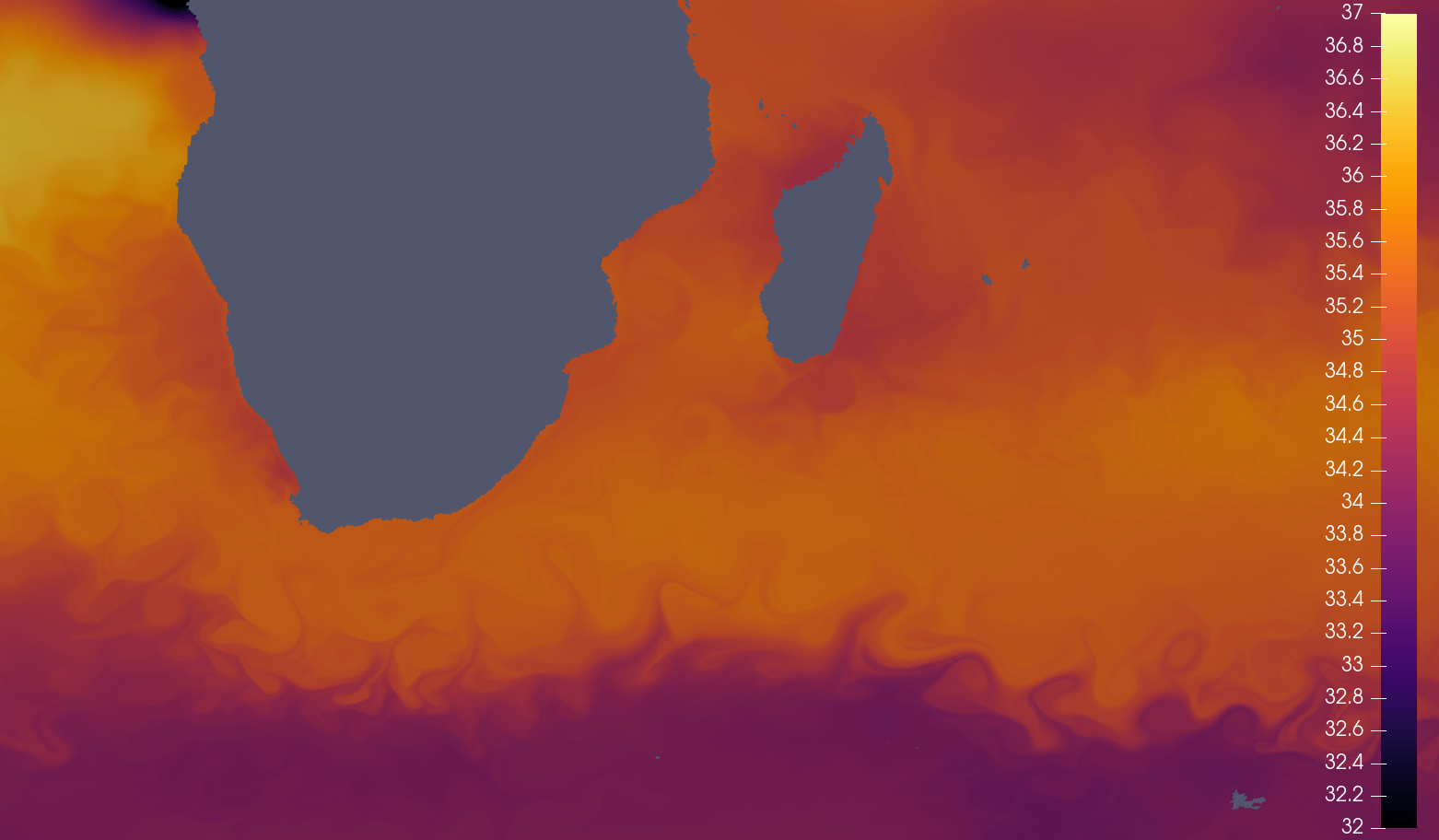}
    \caption{Figures showing average over time of kinetic energy for the reference solution (top left) with surface level temperature (top right), salinity (bottom right), and kinetic energy (bottom left) for an ocean state using AOT. The values for kinetic energy, temperature, and salinity are given with units $m^2s^{-2}$, degrees Celsius, and PSU (Practical Salinity Unit), respectively.}
    \label{fig:4_Piece}
\end{figure}

\begin{figure}
\centering
\includegraphics[width=12cm,height=8cm]{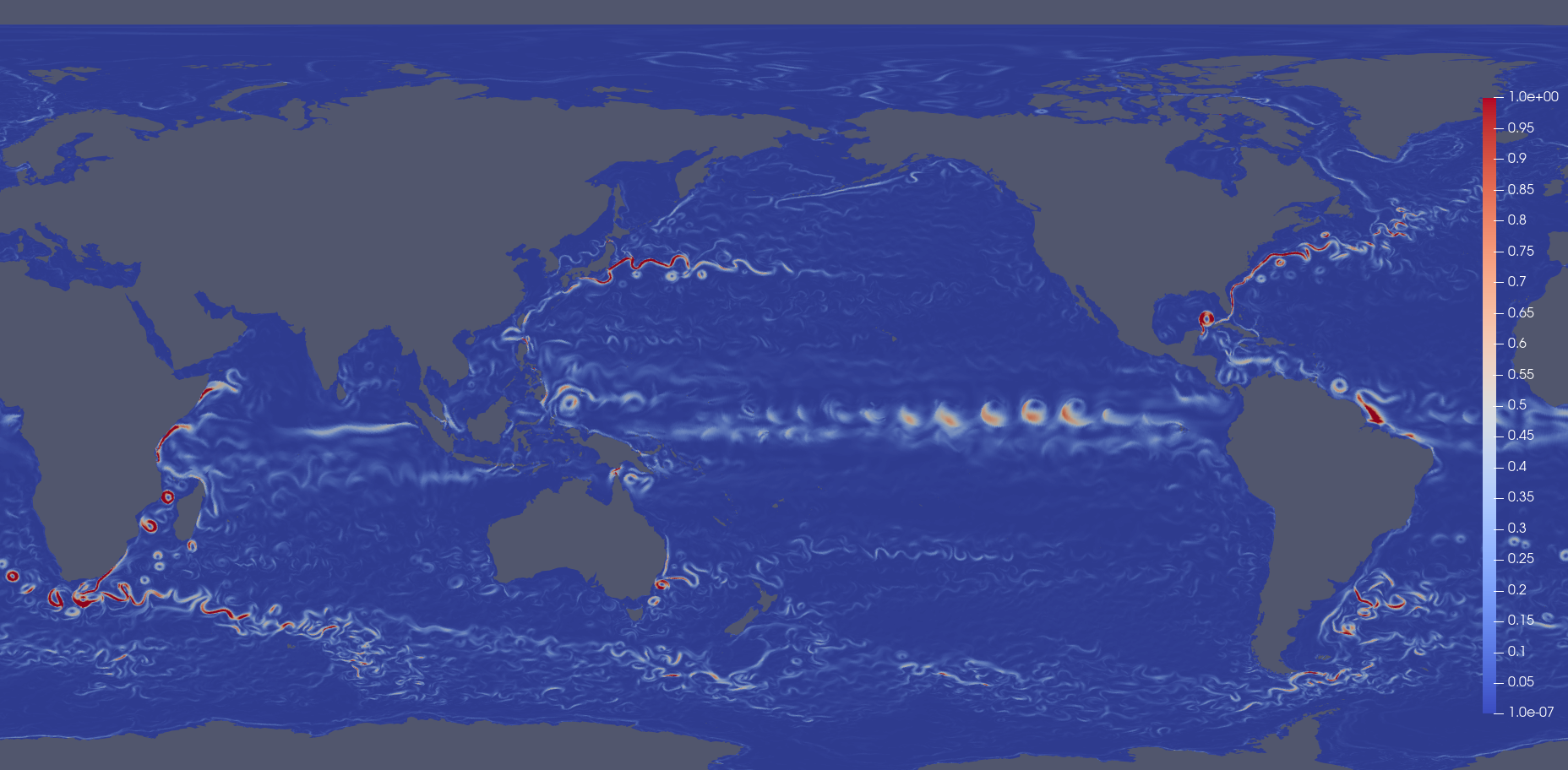}
        \caption{Kinetic energy along the surface of the ocean for simulated solution with optimal parameters after 30 days of simulation. Optimal parameters can be found in \Cref{table: trials}. This figure shows a zoomed out version of the bottom left panel of \Cref{fig:4_Piece} with a rescaled colorbar.}
        \label{full ocean picture}
\end{figure}

\section{Background}\label{background}
In this section we detail the relevant background information related to data assimilation and the ocean model.
We discuss the governing equations for MPAS-Ocean, the ocean model we utilize, the equations for the AOT algorithm applied to this system, and additional details about MPAS-Ocean itself.


In this work all simulations were conducted using the Model for Prediction Across Scale-Ocean (MPAS-Ocean) developed at Los Alamos National Laboratory (LANL).
It is the ocean component of the Department of Energy's Energy Exascale Earth System Model (E3SM) for use in long-term climate simulations. 
MPAS-Ocean utilizes an unstructured mesh to perform ocean simulations with variable spatial resolution. 
This is typically used to conserve computational power by utilizing a high resolution mesh only over areas of particular interest for detailed regional studies, or to resolve smaller-scale features of the flow. The unstructured mesh is initialized using a mesh generation framework known as compass (Configuration Of Model for Prediction Across Scales Setups) which adds bathymetry and initial temperature and salinity data to a layered ocean model. The mesh itself is generated within compass using a Delaunay-based unstructured mesh generation tool called JIGSAW. For an in-depth description of the algorithm used in the mesh generation see e.g. \cite{Engwirda_2018,Engwirda_2016, Engwirda_2015, Engwirda_Ivers_2016, Engwirda_2014}. 
For a more in depth introduction to MPAS-Ocean see, e.g., \cite{mpas-1,Petersen_ea15om,mpas-3}.

\subsection{Governing Equations}
The governing equations for MPAS-Ocean, as given in \cite{mpasoceanuser}, are:
\begin{subequations}
\begin{align}
    \pdv{\bu}{t} + \eta \bk\times \bu + w \pdv{\bu}{z} = -\frac{1}{\rho_0}\grad p - \frac{\rho g}{\rho_0}\grad z^{mid} - \grad K + D_h^u + D_v^u + F^u,\label{momentum}\\
    \pdv{h}{t} + \grad\cdot(h\overline{\bu}^z) + w\rvert_{z = s^{top}} - w\rvert_{z = s^{bot}} = 0,\label{thickness}\\
    \pdv{~}{t} h\overline{\phi}^z + \grad \cdot(h\overline{\phi \bu}^z) + \phi w\rvert_{z = s^{top}} - \phi w\rvert_{z = s^{bot}} = D_h^\phi + D_v^\phi + F^\phi,\label{tracers}\\
    p(x,y,z) = p^s(x,y) + \int_z^{z^s} \rho g dz', \label{hydrostatic}\\
    \rho = f_{eos}(\Theta, S, p).\label{eos}
\end{align}
\end{subequations}
The equations above are the standard expression for the primitive equations. \Cref{momentum} is the momentum equation, \Cref{thickness} is the equation for vertical layer thickness, \Cref{tracers} is the equation for the tracers (salinity and temperature), \Cref{hydrostatic} is the hydrostatic condition, and \Cref{eos} is an equation of state. The variable and parameter values are described in \Cref{table:symbols}.

\begin{table}
\centering
\begin{tabular}{lll}
\hline
symbol               & name                          & notes                                 \\ \hline
%
$D^u_h, D_v^u$       & momentum diffusion terms      & $h$ horizontal, $v$ vertical              \\
$D_h^\phi, D_v^\phi$ & tracer diffusion terms        & horizontal and vertical               \\
$f$                  & Coriolis parameter            &                                       \\
$f_{eos}$            & equation of state             &                                       \\
$F^u$                & momentum forcing              &                                       \\
$F^\phi$             & tracer forcing                &                                       \\
$g$                  & gravitational acceleration    &                                       \\
$\delta$                  & layer thickness               &                                       \\
$\bk$                & vertical unit vector          &                                       \\
$K$                  & kinetic energy                & $K = |\bu|^2/2$                       \\
$p$                  & pressure                      &                                       \\
$p^s$                & surface pressure              &                                       \\
$s^{bot}$            & z-location of bottom of layer &                                       \\
$s^{top}$            & z-location of top of layer    &                                       \\
$S$                  & salinity                      & a tracer $\phi$                                      \\
$t$                  & time                          &                                       \\
$\bu$                & horizontal velocity           &                                       \\
$w$                  & vertical transport            &                                       \\
$z$                  & vertical coordinate           &                                       \\
$z^{mid}$            & z-location of middle of layer &                                       \\
$z^s$                & z-location of sea surface     &                                       \\ 
$\omega$             & relative vorticity            & $\omega = \bk \cdot(\grad \times \bu)$\\
$\eta$               & absolute vorticity            & $\eta = \omega + f$                   \\
$\Theta$             & potential temperature         & a tracer $\phi$                       \\
$\kappa_h, \kappa_v$ & tracer diffusion              & horizontal and vertical               \\
$\rho$               & density                       &                                       \\
$\rho_0$             & reference density             &                                       \\
$\phi$               & generic tracer                & e.g. $\Theta, S$ \\ \hline                     \label{reference}
\end{tabular}
\caption{Variables utilized in MPAS-Ocean governing equations.\label{table:symbols}}
\end{table}

\subsection{AOT Algorithm Equations}
The AOT algorithm is implemented for this system as follows:
\begin{subequations}
\begin{align}
    \pdv{\bu}{t} + \eta \bk\times \bu + w \pdv{\bv}{z} = -\frac{1}{\rho_0}\grad p - \frac{\rho g}{\rho_0}\grad z^{mid} - \grad K + D_h^u + D_v^u + F^u\label{momentum v}\\
     + \mu I_\delta(E_{dt_{obs}}\bu_{ref} - \bu),\nonumber\\
    \pdv{h}{t} + \grad\cdot(h\overline{\bu}^z) + w\rvert_{z = s^{top}} - w\rvert_{z = s^{bot}} = 0,\\
    \pdv{~}{t} h\overline{\phi}^z + \grad \cdot(h\overline{\phi \bu}^z) + \phi w\rvert_{z = s^{top}} - \phi w\rvert_{z = s^{bot}} = D_h^\phi + D_v^\phi + F^\phi \label{tracers v}\\
    + h\mu I_\delta(\overline{E_{dt_{obs}}\phi_{ref}}^z - \overline{\phi}^z)\nonumber, \\
    p(x,y,z) = p^s(x,y) + \int_z^{z^s} \rho g dz',\\
    \rho = f_{eos}(\Theta, S, p). 
\end{align}
\end{subequations}

The variable and parameter values are listed in \Cref{table:aot_variable} as well as explained briefly below. Note that the equations here are the same as the governing equations, except in this case we have additional feedback control terms in \Cref{momentum v} and \Cref{tracers v}. 
In \Cref{momentum v} this additional term is given by $\mu I_\delta(E_{dt_{obs}}(\bu_{ref}) - \bu)$, where $\mu$ is a feedback control parameter, $I_\delta$ is a linear interpolant with spatial resolution $\delta$, $E_{dt_{obs}}$ is piece-wise linear interpolation in time with time-scale $dt_{obs}$, and $\bu_{ref}$ is the value of $\bu$ from our reference solution obtained from \Cref{momentum} that is observed with spatial resolution $\delta$ and time resolution $dt_{obs}$. 
We note that since $I_\delta$ is a linear operator, this term can be equivalently rewritten as $\mu \left(I_\delta(E_{dt_{obs}}(\bu_{ref})) - I_\delta(\bu)\right)$, as in \cite{Carlson_VanRoekel_Petersen_Godinez_Larios_2021}. 
Notice also that we utilize $E_{dt_{obs}}(\bu_{ref})$, which is an approximation of $\bu_{ref}$ created by linearly interpolating between observed states at different times.
The usage of such an approximation is a necessary one, due to the infeasibility of storing high resolution ocean data continuously in time.

It was proved in \cite{Pei_2019} that sufficiently dense observations of the velocity and tracer fields were enough to obtain convergence to the true solution for the primitive equations of the ocean (assuming diffusion in all tracer fields), but observing only the velocity poses analytical problems to the convergence proof. 
Thus, in our setting, it is expected that salinity and temperature measurements will need to be assimilated to achieve convergence of this model to the true solution, and we investigate this below. 
To examine the effect that these tracer quantities have on convergence, we incorporate an additional feedback control term into \Cref{tracers v}.
We note that the feedback control parameter, $\mu>0$, is the same for both \Cref{momentum v,tracers v}.
We found that for simulations the optimal choice of $\mu$ for \Cref{tracers v} appeared to be the same as the optimal choice in \Cref{momentum v} so they are denoted by the same symbol $\mu$.
Note however that in many of our simulations below, we do not assimilate tracer observational data, which can be seen as formally setting $\mu = 0$ only in \Cref{tracers v}.

\begin{table}[]
\centering
\begin{tabular}{lll}
\hline
symbol               & name                          & notes                                 \\ \hline
%
$I_\delta$                & interpolation operator        & spatial length scale $\delta$              \\
$E_{dt_{obs}}$                & interpolation operator        & time scale $dt_{obs}$              \\
$\delta$                  & spatial length scale          & measured in cell lengths              \\
$dt_{obs}$                  & time scale          & time between observations \\
$\mu$  & feedback control parameters   & \\
$\bu_{ref}, \phi_{ref}$& reference solution observables          & observed with spatial resolution $\delta$ \\ \hline \label{simulated}
\end{tabular}
\caption{Variables utilized by the AOT algorithm formulation applied to the governing equations.\label{table:aot_variable}}
\end{table}


\section{Methods}\label{methods}

In this section we outline the experimental setup that we employed to test the AOT algorithm in the MPAS-Ocean system. Additionally we detail the numerical schemes and algorithms used in the simulation of the governing equations and for the generation of observational data. 
We also discuss several validation tests that were utilized to ensure that the AOT algorithm was implemented correctly into MPAS-Ocean.


\subsection{Experimental Setup}
In order to test the AOT algorithm we used an ``identical twin'' experimental design. This is a standard experiment used to test methods of data assimilation, involving running two separate simulations. The first simulation is used to produce a reference solution, which is used both to generate observational data and as the ``true'' solution we want to recover using data assimilation. Next, we start the twin simulation, a second simulation that we begin from a different set of initial data. 
The observational data is assimilated into the twin simulation, and  the error between both simulations is calculated. 
The phrase identical twin here refers to the underlying physics in the simulations; namely, both the reference and twin simulation utilize the same underlying model with consistent physics between them.
We note that in several validation tests the physics were adjusted (as detailed below), which would classify those tests as ``fraternal twin'' experiments instead. 
For a detailed discussion of twin experiments see, e.g., \cite{yu2019twin}.

For our experiments we first initialized the reference solution using high resolution initial data which was ``spun-up'' until the kinetic energy was statistically stabilized. 
This data was then run for another two simulation years to ensure that the kinetic energy was approximately constant for the system. 
The end result of these simulations was then utilized as the initial state for the reference solution.
This reference solution was simulated for an additional month in order to generate observational quantities.
The final state of this simulation was recorded and utilized as the initial state for the twin simulation.
We note that after thirty days have passed, the end state of the reference solution is sufficiently decorreleated from the initial state to ensure that the root-mean-squared error of the kinetic energy between the simulated solution and the reference solution is large (e.g. $5e-2$)\footnote{We use the standard notation $Me-n$ to mean $M\times10^{-n}$.} during the initial startup of the simulated solution.
Moreover the solutions are decorrelated to the point where the error remains large for all times when no observational data is assimilated.

We note that we utilized the difference of pointwise kinetic energy as a metric for the error as it is a standard metric in oceanography used to characterize simulations. In addition, the kinetic energy is a scalar and so does not depend on the frame of reference, as the velocity components do; which have meridional and zonal components when measured in cell centers or normal and tangential components when measured on cell edges. In low resolution tests, we found that the differences between in error between simulations was qualitatively the same when using the kinetic energy and the velocity fields directly.


\subsection{Numerical Scheme}

MPAS-Ocean was developed with the option to use different time-stepping schemes to advance the governing equations forward in time. 
The methods we used in our simulations were a barotropic-baroclinic splitting explicit scheme (SE), as well as a forward Euler scheme that we implemented for our validation tests. 
MPAS-Ocean additionally has a fourth order Runge-Kutta (RK4) scheme, however we note that, as pointed out in \cite{Olson_Titi_2008_TCFD}, the AOT algorithm cannot be faithfully implemented with higher-order Runge-Kutta methods, or implicit methods, due to the fact that the feedback control term, being based on data, cannot be evaluated at arbitrary arguments. 
We note that the forward Euler scheme is applied specifically to the momentum equation. We found that the inclusion of terms in the tracer and vertical layer thickness equations introduced instabilities, making the forward Euler scheme unsuitable for the high resolution tests using the full ocean model. 

Details of the specific time-stepping schemes used in MPAS-Ocean can be found in \cite{mpas-1}. We note that the forward Euler scheme was implemented as a modification of the current RK4 scheme. 
The assimilated equations were simulated using the same numerical scheme used to generate the initial data, with the exception that the feedback control term was implemented as an explicit term in both the momentum and tracer equations, except for certain validation tests which we discuss below.


Additionally we note that our simulations were all conducted on the Cori computing cluster at the National Energy Research Scientific Computing Center (NERSC). 
We utilized Cori's Haswell system partition which utilizes Intel Xeon processors (E5-2698 v3) running at a clock rate of 2.3 GHz. Numerical tests were performed using 512 cores for low resolution and 4096 cores for high resolution simulations of the ocean.



\subsection{Observational Data and Interpolant}
The observational data was gathered by saving the meridional and zonal velocities and the salinity and temperature tracers as measured in the cell centers at 6-hour intervals. Low resolution simulations showed that the difference in convergence between simulations using observations gathered daily and hourly were negligible. Due to this, the choice of 6-hour time intervals was made as a compromise to conserve storage space while still utilizing a significant amount of observational data.
We note that for ultra-low resolution simulations, we ran tests with ocean data saved from each time-step. The results of these simulations were qualitatively the same as those featuring a linear interpolation of observational data for the high resolution simulations, so they have been omitted from this work.

The linear interpolant we used for $I_\delta$ was a modified flood fill algorithm given in \Cref{alg:flood} with an additional smoothing given in \Cref{alg:smoother}. First, we used a masking array to determine the locations of observational data. Using this mask, we initialized the observed cells with the meridional and zonal velocities as measured in the cell centers. Once the observation cells are initialized we perform a flood fill algorithm. That is, we propagate all of the observational data from the unmasked cells to all neighboring cells. This propagation continues from all newly initialized cells to all of its neighbors. This continues until all cells in the domain are initialized. In the case of a conflict, where multiple cells try to propagate their information to a shared neighbor, the data in the new cell is averaged between all neighboring cells that have been initialized in the previous iteration. 
After the observational data has been initialized using this algorithm the data in all of the unobserved cells is then averaged with all of the neighboring cells using \Cref{alg:smoother}. This averaging is done to ensure a smoother transition between the piecewise-constant regions around the observed cells.

\vfill

\phantom{a}

\begin{algorithm}[H]
\KwData{Ocean observation mask: $Obs(i) = \begin{cases}
    1 & \text{ cell $i$ is observed,}\\
    0 & \text{ cell $i$ is not observed}
\end{cases}.$}
\KwData{Current ocean data on cells, $O_i$, given over $N$ total cells.}
\KwData{Number of smoothing iterations, $smoothing$, betweeen $0$ and $\delta$.}
\KwResult{Interpolated ocean data $\tilde{O}_i$}
\For{$i = 1:N$}{
    \%\% Initialize variable tracking distance of individual cells from observed cells.\\
    $cell\_depth(i) = \begin{cases}0 & Obs(i) = 1,\\ \infty & Obs(i) = 0 \end{cases}$\;
    \%\% Initialize interpolated ocean with data from observed cells\\
    $\tilde{O}_i = Obs(i)\cdot O_i$\;
    \%\% Initialize variable that maintains a count for each cell used for averaging.\\
    $count(i) = 0$\;
}
\%\% Loop over cells targeting a specific distance from observed cells ($depth$).\\
\For{$depth = 0$ to $\delta$}{
    \For{$i = 1:N$}{
        \eIf{$cell\_depth(i) = depth$}{
            \%\% Extract list of neighboring cells\\
            $nearby\_cells = Neighbors(i)$\;
            \%\% Propogate information to neighboring cells\\
            \For{$cell$ in $nearby\_cells$}{
                \eIf{$cell\_depth(cell) > depth$}
                    {
                    \%\% Set correct depth\\
                    $cell\_depth(cell) = \min(cell\_depth(cell), depth + 1)$\;
                    \%\% Update $\tilde{O}$ with running average.\\
                    $\tilde{O}(cell) = average(\tilde{O}(cell), \tilde{O}(i), count)$\;
                    \%\% Update count for running average.\\
                    $count(cell) = count(cell) + 1$\;
                    }
                    {
                    \%\% Do nothing. Cell has been flooded by closer cells already. 
                    }
            }
        }
        { \%\% Do nothing. Cell is not at the targeted distance from observed cells.
        }
        
    }
}
\For{$j = 1: smoothing$}{
\%\% Apply smoothing algorithm.\\
    $\tilde{O} = smoothing(\tilde{O})$\;
}
\caption{Flood Fill Interpolant}
\label{alg:flood}
\end{algorithm}

\begin{algorithm}
\KwData{Interpolated observational data to be smoothed: $\tilde{O}$.}
\KwData{Ocean observation mask: $Obs(i) = \begin{cases}
    1 & \text{ cell $i$ is observed,}\\
    0 & \text{ cell $i$ is not observed}
\end{cases}.$}    
\%\% Save data to new variable to prevent race conditions.\\
$\tilde{O}\_temp = \tilde{O}$\;
\%\% Loop over all $N$ cells in the domain.\\
    \For{$i=1:N$}{
        \eIf{$Obs(i) \neq 1$}{
        \%\% Find neighboring cells of current cell\\
            $nearby\_cells = Neighbors(i)$\;
            $total = \tilde{O}\_temp(i)$\;
            $count = 1$\;
            \%\% Sum over all neighboring cells\\
            \For{$cells$ in $nearby\_cells$}{
                $total =  total + \tilde{O}\_temp(cells)$\;
                $count = count + 1$\;
            }
           $\tilde{O}(i) = total/count$\;
        }{
        \%\% Do nothing, as these cells are directly observed.
        }
            
    }
\caption{Smoothing Operator}
\label{alg:smoother}
\end{algorithm}

\begin{figure} 
        \includegraphics[width=4.5cm,height=3cm]{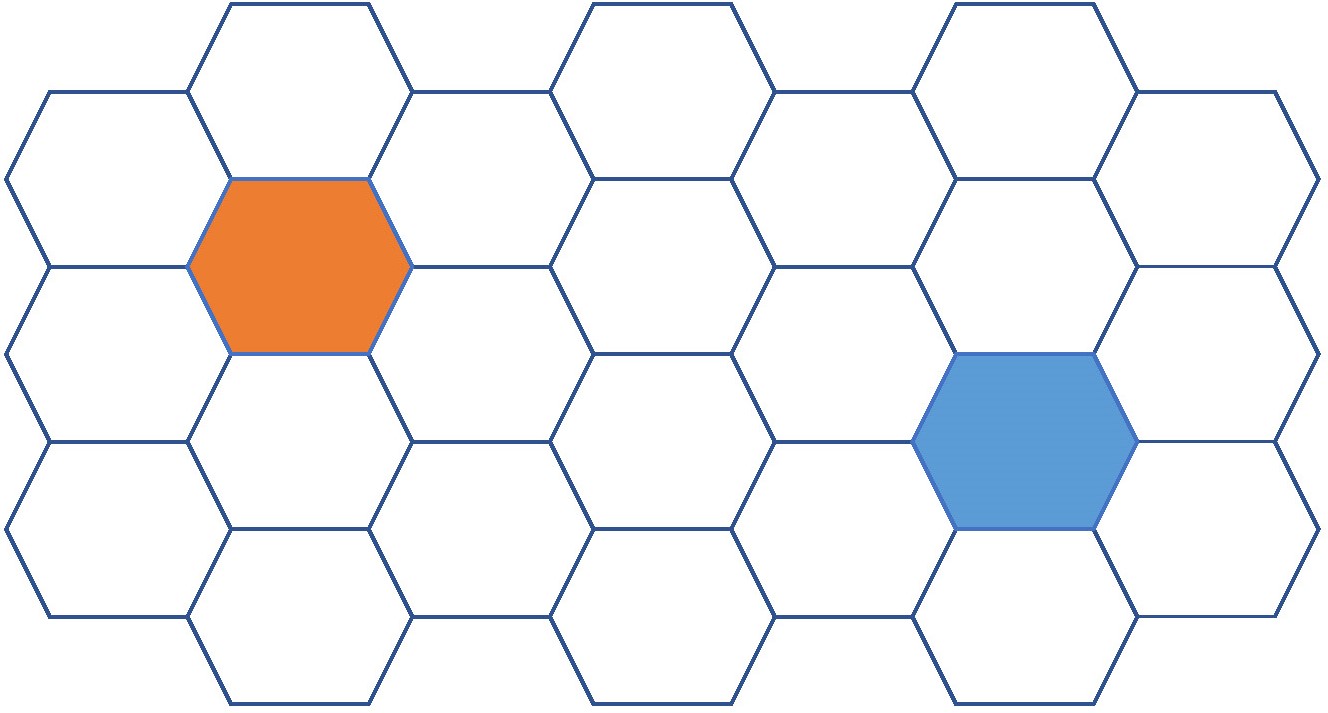}
        \includegraphics[width=4.5cm,height=3cm]{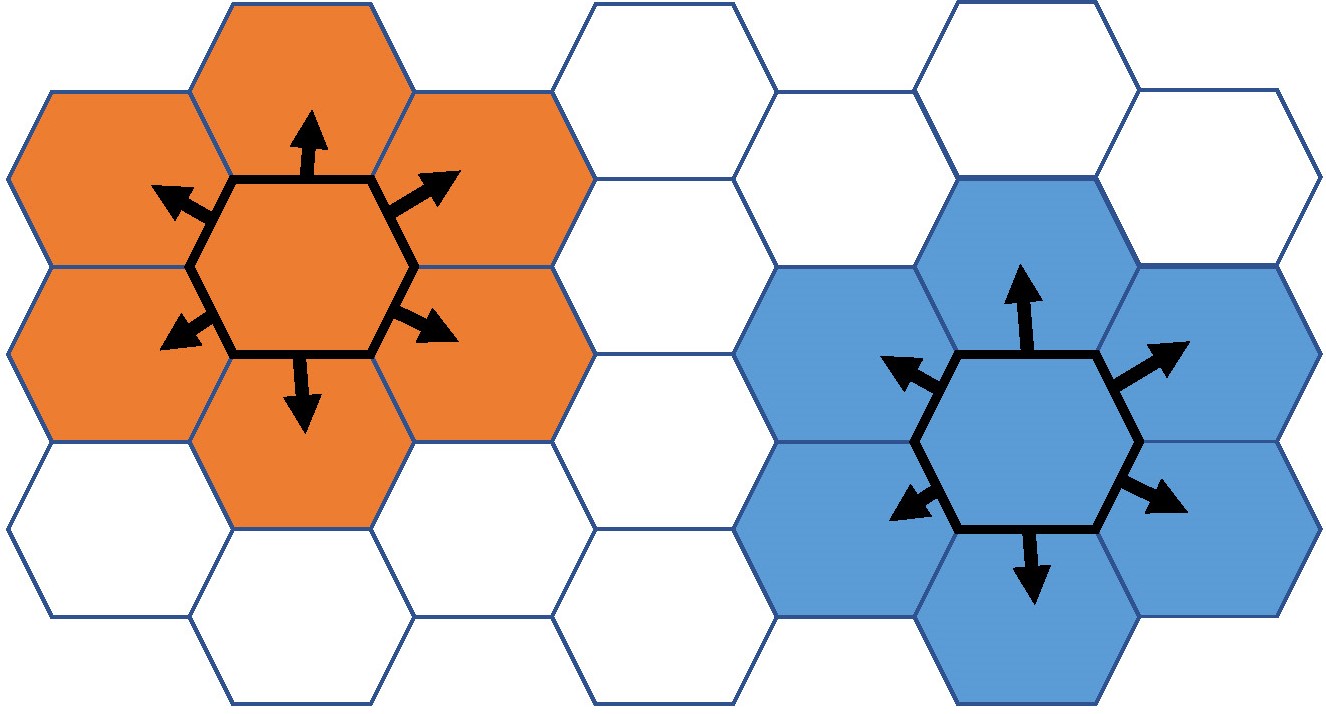}
        \includegraphics[width=4.5cm,height=3cm]{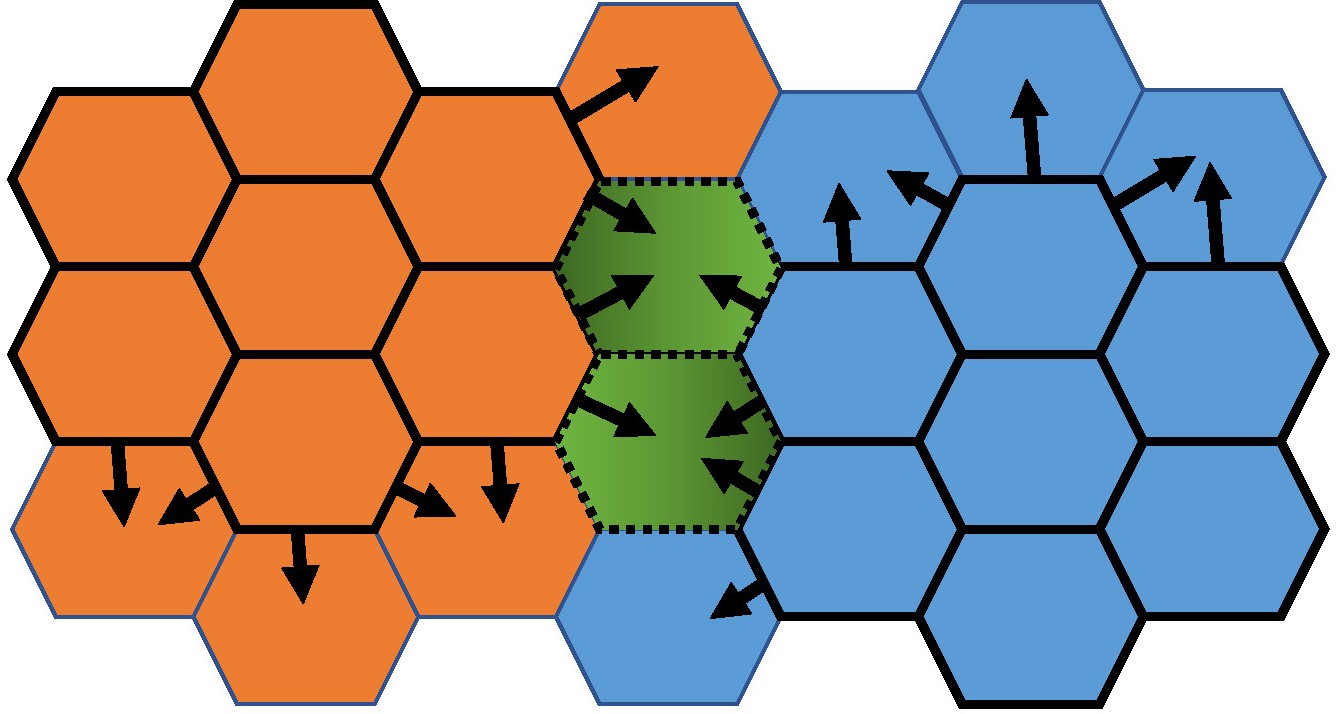}
\caption{Diagram of flood fill interpolation algorithm. Observation data is first initialized with data (given by colorings) from observed cells and zero on unobserved cells (left), then the algorithm loops based on distance from observed cells (middle and right), where information from active cells is propagated to neighboring cells via arrows. In case of conflicts cell data is averaged over all cells attempting to propagate to it as indicated via a green coloring in the central cells in the rightmost image.}
\end{figure}
\subsection{Validation Tests}
To validate the implementation of the AOT algorithm, validation tests were performed on low resolution simulations. These tests were conducted to illustrate that the algorithm was implemented correctly despite the fact that the convergence to machine precision was not achieved using observational data of arbitrary spatial and temporal distributions. The validation tests were conducted to test the implementation for different time-stepping lengths, convergence to a static ocean state, and implicit/explicit variations of the AOT algorithm. The validation tests featured in this section were performed on a low resolution ocean mesh, using 512 processor cores in simulations.

     \begin{figure}
        \centering
        \includegraphics[width=12cm,height=8cm]{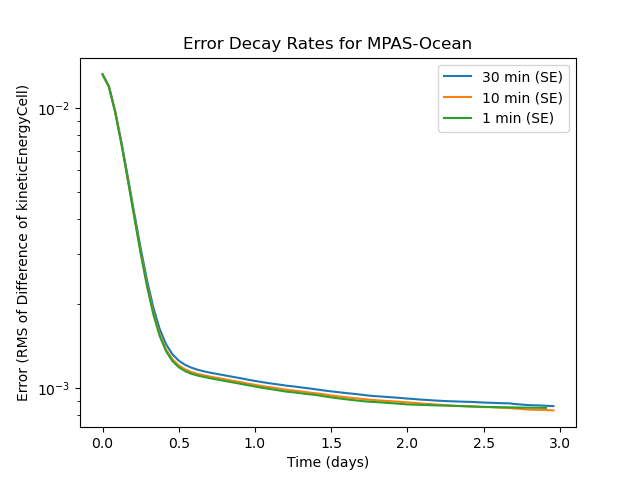}
         \caption{Convergence rates for varying $dt$ choices. The "SE" (Split Explicit) in the legend refers to the barotropic baroclinic splitting scheme used in the time integration.}
         \label{fig:valid dt}
     \end{figure}

     \begin{figure}
        \centering
      \includegraphics[width=12cm,height=8cm]{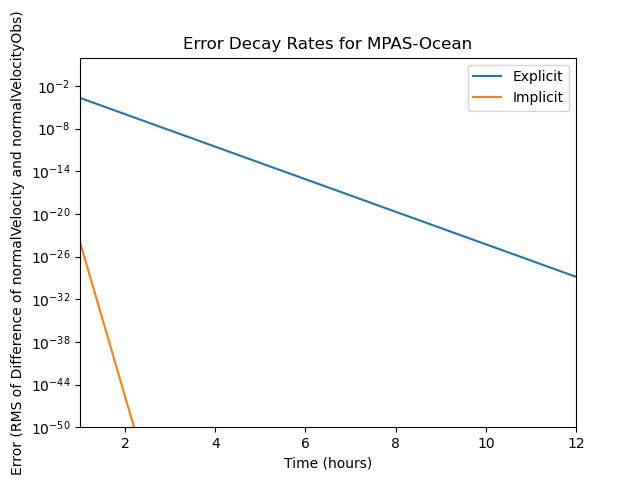}
           \caption{Convergence to static ocean state for an implicit and explicit implementations of the AOT term.  Note: The output times for MPAS-Ocean are much larger than the time-step, and hence the above plot only shows the error after the first hour of simulation.}
         \label{fig:valid static}
     \end{figure}

\begin{table}[]
\centering
\begin{tabular}{llll}
\hline Trial & Explicit & Implicit\\
\hline
No Dynamics & EPS & EPS\\
Full Dynamics & 3.107e-4 & 2.515e-4\\
Coriolis Only & 1.090e-4 & 1.781e-4\\
Explicit Bottom Drag Only & EPS & EPS\\
Surface Stress Only & 1.488e-8 & 3.721e-12\\
Topographic Wave Drag Only & EPS & EPS\\
Horizontal Mixing Only & 7.672e-6 & 1.919e-9\\
Vertical Mixing Only & 7.535e-4 & 1.083e-3\\
Vertical Advection Only & 9.695e-8 & 3.264e-11\\
Pressure Gradient Only & 1.379e-4 & 3.449e-8\\
\hline 
\end{tabular}
\caption{Convergence levels for implicit and explicit implementations of the AOT feedback control term into the forward Euler numerical scheme for the momentum equation. Error is measured as the RMS of the difference between the normal velocity and the observed normal velocity. All implicit simulations used $\mu = 10$ and all explicit simulations used $\mu = 1e-3$. Error level was achieved within 3 hours of simulation and stayed approximately constant. EPS for certain explicit trials indicates convergence to machine precision with the error returning a NAN (not a number) error when calculated at the end of the 24 hour simulation due to lack of precision.} \label{table: schemes}
\end{table}

In order to ensure the efficacy of these validation tests we adjusted the implementation of the AOT algorithm to use edge-based measurement data. Normally the observational data is recorded at the cell centers for all of the cells on the unstructured mesh utilized by MPAS-Ocean. This presents a potential source of error, as the governing equations are evolved using quantities evaluated on the cell edges, which requires the observational data to interpolated from cell centers to cell edges. This interpolation breaks down at the poles, which is why we implemented a latitude filter to only use observational data below $70^\circ$ North. For the validation tests we instead use edge-based measurements directly to avoid any error from interpolation. 


Additionally, in all of our validation tests we set the observational length scale of ``$\delta = 0$.'' This means that we are utilizing the full solution $\bu_{ref}$ as observational data. This is useful for validating the methods, however for this case the algorithm does simplify to the more classical data assimilation known as nudging, as there is no spatial interpolation of the observational data.

\begin{figure}
    \centering
        \includegraphics[width=12cm,height=8cm]{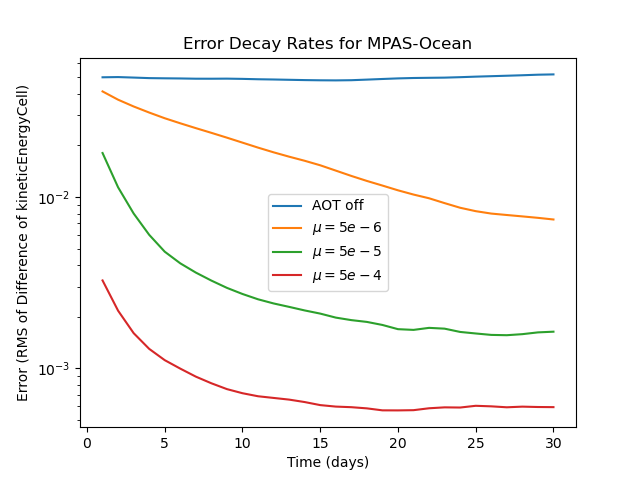}
\caption{Error corresponding to varied $\mu$ values over time.}
    \label{fig:mu}
\end{figure}

The results of our validation tests can be seen in \Cref{fig:valid dt,,fig:valid static,,table: schemes}. The first validation test we conducted was an alteration of the time-stepping parameter $dt$. 
This was done to ensure that the results we obtained were consistent for different time-steps that shared the same $\mu$ value. 
Due to the explicit implementation of the feedback control term, there is a CFL-type condition restricting the choice of $\mu$ parameter; roughly, $\mu  \lesssim \frac{2}{dt}$. 
This allows smaller time-steps to utilize a smaller $\mu$ value for potentially faster convergence.
In our simulations we saw that, as expected, utilizing the same $\mu$ value across a selection of $dt$ values resulted in similar error rates (see \Cref{fig:valid dt}). 
Moreover, we also noted that if we set $\mu = \frac{\mu_0}{dt}$ with $\mu_0$ a constant, then the scheme remained stable with better error values, as expected.
This was done to ensure the CFL condition is satisfied while increasing the value of $\mu$ for smaller time-steps.
The results of these simulations are not pictured, but are qualitatively similar to \Cref{fig:mu}.
While decreasing the time-step may be somewhat beneficial to accuracy, it means that one is spending significantly more computational resources simulating the equations. 
Additionally, we found that we were unable to achieve convergence to machine precision using this method, even with a time-step far smaller than needed for low resolution simulations.

Another validation test we performed was to test convergence to a static ocean state (i.e., a time-independent state). In this case we used a static ocean state as our reference solution $\bu_{ref}$ and attempted to obtain convergence to it using the AOT algorithm. For these tests, we turned the underlying ocean physics off, as the static ocean state is arbitrarily chosen and not in equilibrium. In this case the AOT equations reduce from a PDE to an ODE due to the modifications to the model equations, greatly simplifying the dynamics. The results of this test can be seen in \Cref{fig:valid static}. Note that here we use the forward Euler scheme with either an implicit or explicit implementation of the feedback control term, which we will discuss more in depth later on. 
It is worth mentioning that the forward Euler scheme was implemented for this test instead of the split-explicit scheme used in other tests as the split-explicit scheme does not allow for disabling of individual ocean terms in the way we require. 
Notice in \Cref{fig:valid static} that we observe exponential convergence in time with arbitrary precision. Machine-level precision ($\sim$1e-15) is achieved within one hour of simulation time for the implicit scheme and 8 hours for the explicit scheme.

One additional validation test that we performed was a modification of the implementation of the AOT term. 
We reformulated the algorithm into a predictor-corrector scheme by treating the $\mu I_\delta v$ term implicitly and the $\mu I_\delta E_{dt_{obs}} u$ term explicitly. We note that, since there is no interpolation for this validation test, this scheme is equivalent to a semi-implicit implementation of nudging, a more classical data assimilation algorithm. 
This test was inspired by our previous observation that increasing the value of $\mu$ improves convergence levels of simulations.
By implementing the $\mu I_\delta v$ term implicitly, we could relax stability constraint $\mu  \lesssim \frac{2}{dt}$ and choose a much larger $\mu$ value for a given time-step.  Note that implicit methods for the AOT algorithm were also investigated in \cite{Larios_Rebholz_Zerfas_2018} (but, as a cautionary note, see the discussion in Section 4 of \cite{Olson_Titi_2008_TCFD}).

To see where the CFL constraint $\mu \lesssim \frac{2}{dt}$ arises, consider the ODE
\begin{align}
    v_t = -\mu(v). 
\end{align}
Approximating the time derivative with the explicit Euler scheme we obtain:
\begin{align}
    \frac{v^{n+1} - v^{n}}{dt} = -\mu v^{n},
\end{align}
where the indices refer to the timestep $v$ is evaluated at, i.e. $v^n = v(t_n)$.
This simplifies to the following equation
\begin{align}
    v^{n+1} &= (1-dt\mu)v^{n}\\
    &= (1-dt\mu)^n v^0.
\end{align}
Thus, to ensure that our scheme is stable we require $\abs{1-dt\mu} \leq 1$, which implies $\mu \leq \frac{2}{dt}$, for $dt$ and $\mu$ positive.
Now, if we instead use the implicit Euler update scheme we obtain
\begin{align}
    \frac{v^{n+1} -v^n}{dt} = -\mu v^{n+1}
\end{align}
which simplifies to 
\begin{align}
    v^{n+1} &= (1+dt\mu)^{-1}v^{n}\\
    & = (1 + dt\mu)^{-n}v^0. \label{Implicit CFL}
\end{align}
Now similarly to before, we require $\abs{(1+dt\mu)^{-1}} \leq 1$ for the scheme to remain stable, however note now that this constraint is automatically satisfied for any positive $dt$ and $\mu$.

We simulated the new implicit implementation, comparing it to the original explicit scheme to observe any changes in performance. The results of these simulations can be seen in \Cref{table: schemes}. In addition to testing the implementation, we also investigated the performance of both schemes with different physics enabled. We found that machine-level precision can be achieved for both schemes with certain physical processes disabled. The Coriolis force and vertical mixing terms appeared to be particularly restrictive, as they restrict the error to $1e-4$. The pressure gradient term introduced a similar level of error to these terms for the explicit implementation, but was much less restrictive for the implicit implementation. 

We note that each of the simulations involved in this validation test exhibit the same qualitative behavior in the error plots, where the error very quickly drops to a constant value where it remains for the rest of the simulation. 
We were able to achieve convergence to machine precision for the both schemes, but this was only realized in the case when disabling all of the physical ocean terms except those involving wave drag. 
It appears from these simulations that the Coriolis and vertical mixing terms are especially problematic, as these terms have errors significantly larger than any of the other terms. The only terms that are not problematic for convergence are topographic wave drag and explicit bottom drag.
When the physics of the model are adjusted, the error of implicit and explicit schemes still behave the same qualitatively. 
That is, each simulation exhibited steep decay in the error that leveled off quickly and remained roughly constant for the remainder of the simulation time.
We note that the error level achieved for explicit and implicit methods are close when simulating all of the ocean physics, even though the feedback control parameter $\mu$ is much larger for the implicit schemes.

\section{Computational Results}\label{computational results}

In this section we will outline the computational results we obtained in this study. 
In particular we look at sensitivity of the error to several of the parameters present in the assimilated system. 
Using high resolution simulations we simulated the AOT algorithm as applied to the MPAS-Ocean system. 
We first look at the application of the AOT algorithm applied to the momentum equations exclusively. 
We note that to apply the current theory (see \cite{Pei_2019}) in order to achieve exponential convergence to arbitrary precision for the primitive equations of the ocean, we require assimilation terms for both momentum and tracer quantities, at least in an idealized setting without salinity effects. 
In our study we examine the usefulness of data assimilation applied only to the momentum equation as well as with the inclusion of tracer quantity data assimilation. 

\begin{center}
    \begin{figure}
        \centering
        \includegraphics[width=12cm,height=8cm]{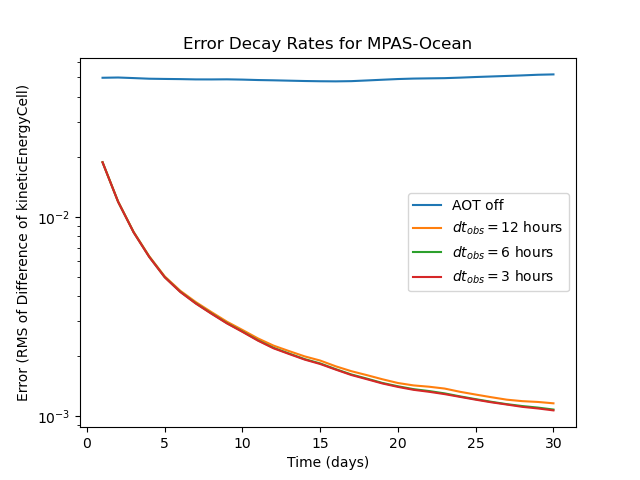}
        \caption{Error corresponding to varied $dt_{obs}$ values over time.}
        \label{dt obs}
    \end{figure}
\end{center}

\begin{center}
    \begin{figure}
        \centering
        \includegraphics[width=12cm,height=8cm]{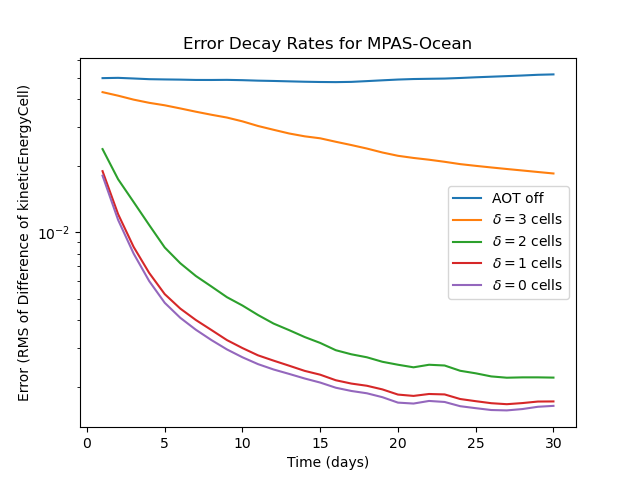}
        \caption{Error corresponding to varied $\delta$ values over time.}
        \label{fig:h}
    \end{figure}
\end{center}

\begin{center}
    \begin{figure}
        \centering
        \includegraphics[width=12cm,height=8cm]{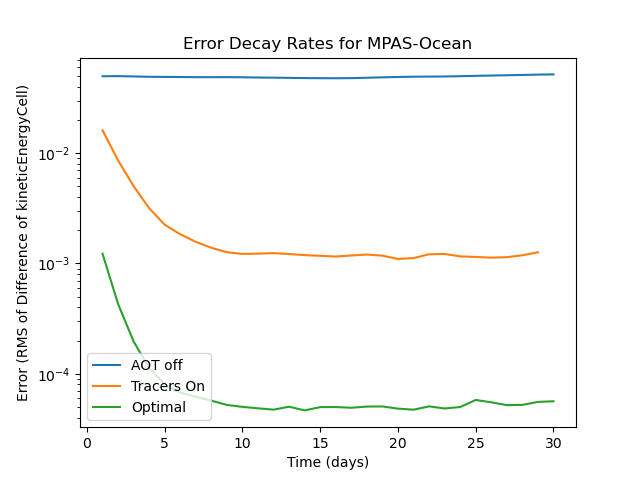}
        \caption{Error corresponding to default and optimized parameter values based on previous tests.}
        \label{fig:optimized}
    \end{figure}
\end{center}


\subsection{Parameter Dependence}

For this study we varied a variety of parameters relevant to the AOT algorithm. These parameters were $\mu$, $dt_{obs}$, and $\delta$.  Here, $\mu$ is the feedback control parameter from the AOT equations, $dt_{obs}$ is the length of time between observational data is recorded, and $\delta$ is the spatial length scale between observer locations (measured in cell lengths). The errors shown in \Cref{table: trials} were calculated as the maximum value of the root mean square of the difference between the kinetic energies calculated between the reference and simulated solutions taken over daily measurements calculated over 30 simulated days. We note that the simulations featured here were conducted on the high resolution ocean mesh, using 4096 cores.

It is important to note that in our above formulation, the AOT algorithm utilizes observational data that is continuous in time (although, see, e.g., \cite{Larios_Pei_Victor_2023_second_best,Foias_Mondaini_Titi_2016,Celik_Olson_2022,Jolly_Martinez_Olson_Titi_2018_blurred_SQG} for the AOT algorithm adapted to discrete times). In our case, the generation and storage of high resolution observational data presents a problem for the application of the AOT algorithm as is. To solve this issue we introduced a new parameter, mentioned above, known as $dt_{obs}$. Instead of recording observational data at each time step, we instead record observational data at each $dt_{obs}$ time interval and performed a linear interpolation in time. In low-resolution tests it was found that the value of $dt_{obs}$ was largely irrelevant to the convergence level attained. This can be seen in \Cref{dt obs} where the plots for 3 hour and 12 hour values of $dt_{obs}$ results in the errors fall roughly on top of one another. 


\begin{table}[]
\centering
\begin{tabular}{llllll}
\hline
Trial               & $\mu$                          & $dt_{obs}$ & $\delta$ & Tracers & Error                                \\ \hline
$\mu$ \#1                      & 1e-4                                    & 03:00                            & 1                      & Off                          & 1.13e-3                    \\
$\mu$ \#2                      & 1e-5                                    & 03:00                            & 1                      & Off                          & 3.60e-3                    \\
$\mu$ \#3                      & 1e-6                                    & 03:00                            & 1                      & Off                          & 5.48e-3                    \\
$dt_{obs}$ \#1                 & 1e-5                                    & 03:00                            & 1                      & Off                          & 1.64e-3                    \\
$dt_{obs}$ \#2                 & 1e-5                                    & 06:00                            & 1                      & Off                          & 3.60e-3                    \\
$dt_{obs}$ \#3                 & 1e-5                                    & 12:00                            & 1                      & Off                          & 1.53e-3                    \\
$\delta$ \#1                & 1e-5                                    & 03:00                            & 0                      & Off                          & 3.60e-3                    \\
$\delta$ \#2                & 1e-5                                    & 03:00                            & 1                      & Off                          & 3.60e-3                    \\
$\delta$ \#3                & 1e-5                                    & 03:00                            & 2                      & Off                          & 3.74e-3                    \\
Tracers  Off                & 1e-5                                    & 03:00                            & 1                      & Off                          & 3.60e-3                    \\
Tracers On                  & 1e-4                                    & 03:00                            & 1                      & On                           & 5.853e-4                    \\
Default                     & 1e-5                                    & 03:00                            & 1                      & Off                          & 3.60e-3                    \\
Optimized                   & 1e-4                                    & 03:00                            & 0                      & On                           & 1.77e-4                    \\ \hline
\end{tabular}
\caption{Convergence levels for various choices of essential data assimilation parameters.}
\label{table: trials}
\end{table}

 \begin{figure}
    \centering    
    \includegraphics[width=.3\textwidth, height=4cm]{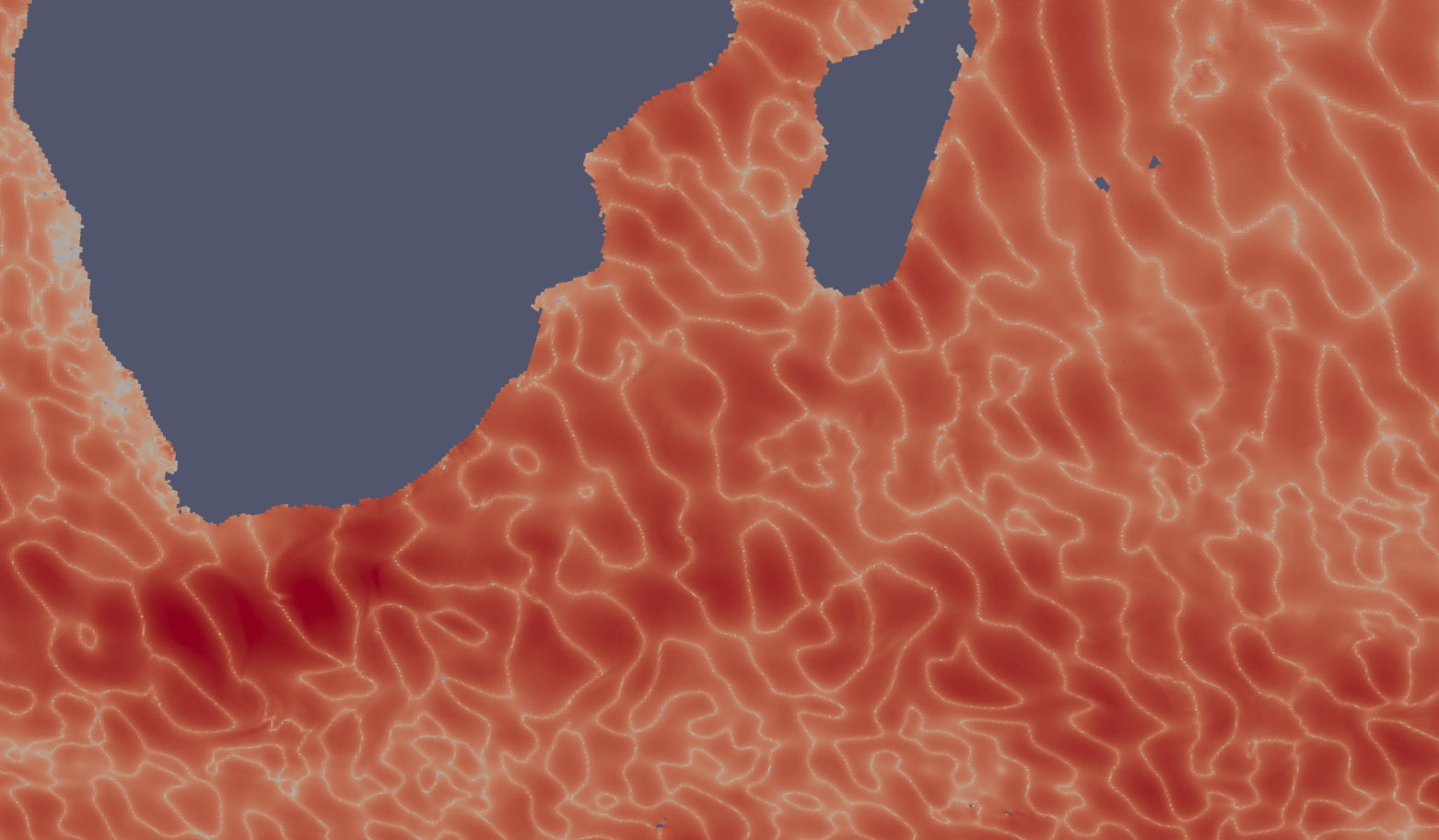}    
    \includegraphics[width=.3\textwidth, height=4cm]{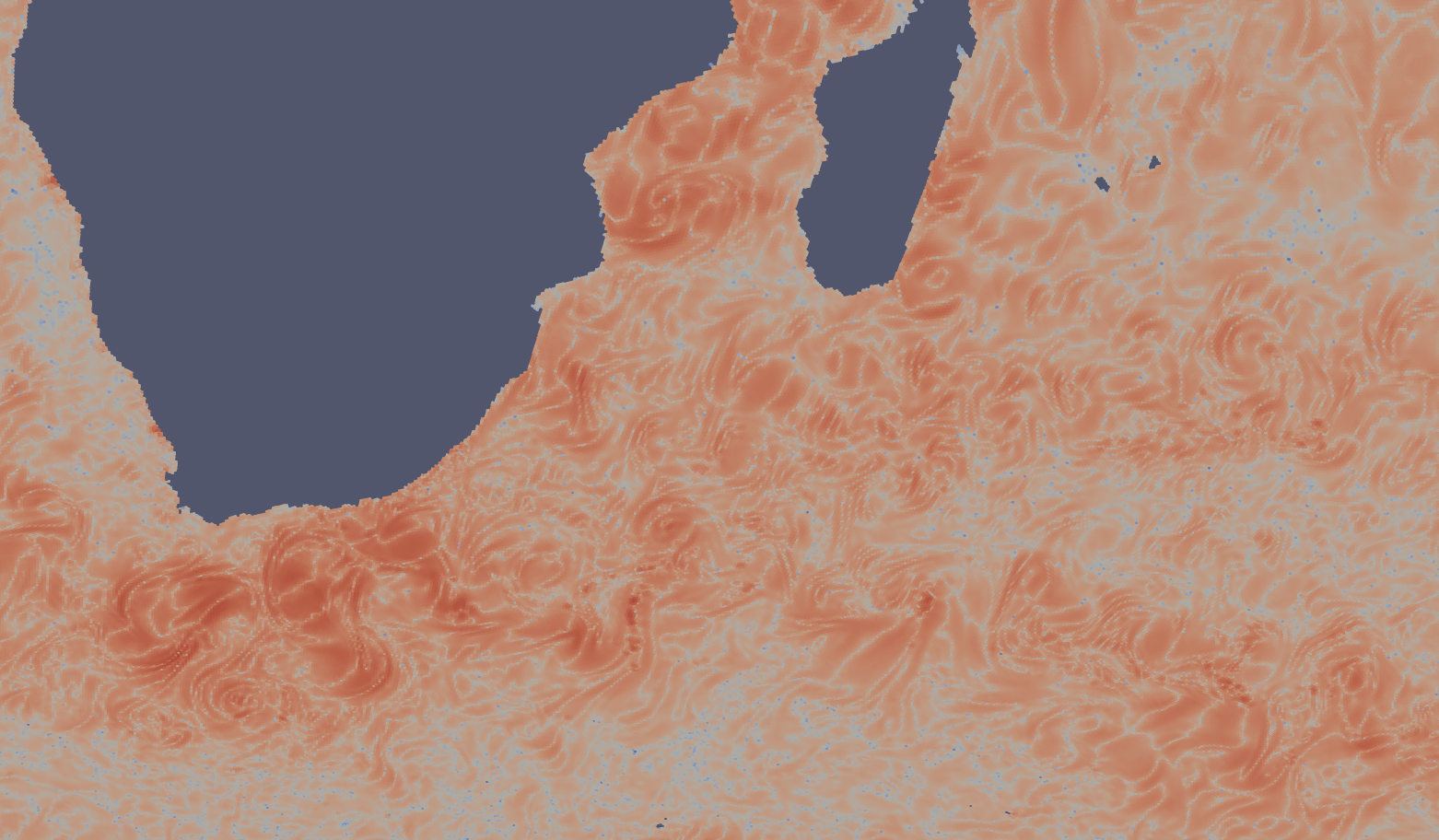}    
    \includegraphics[width=.3\textwidth, height=4cm]{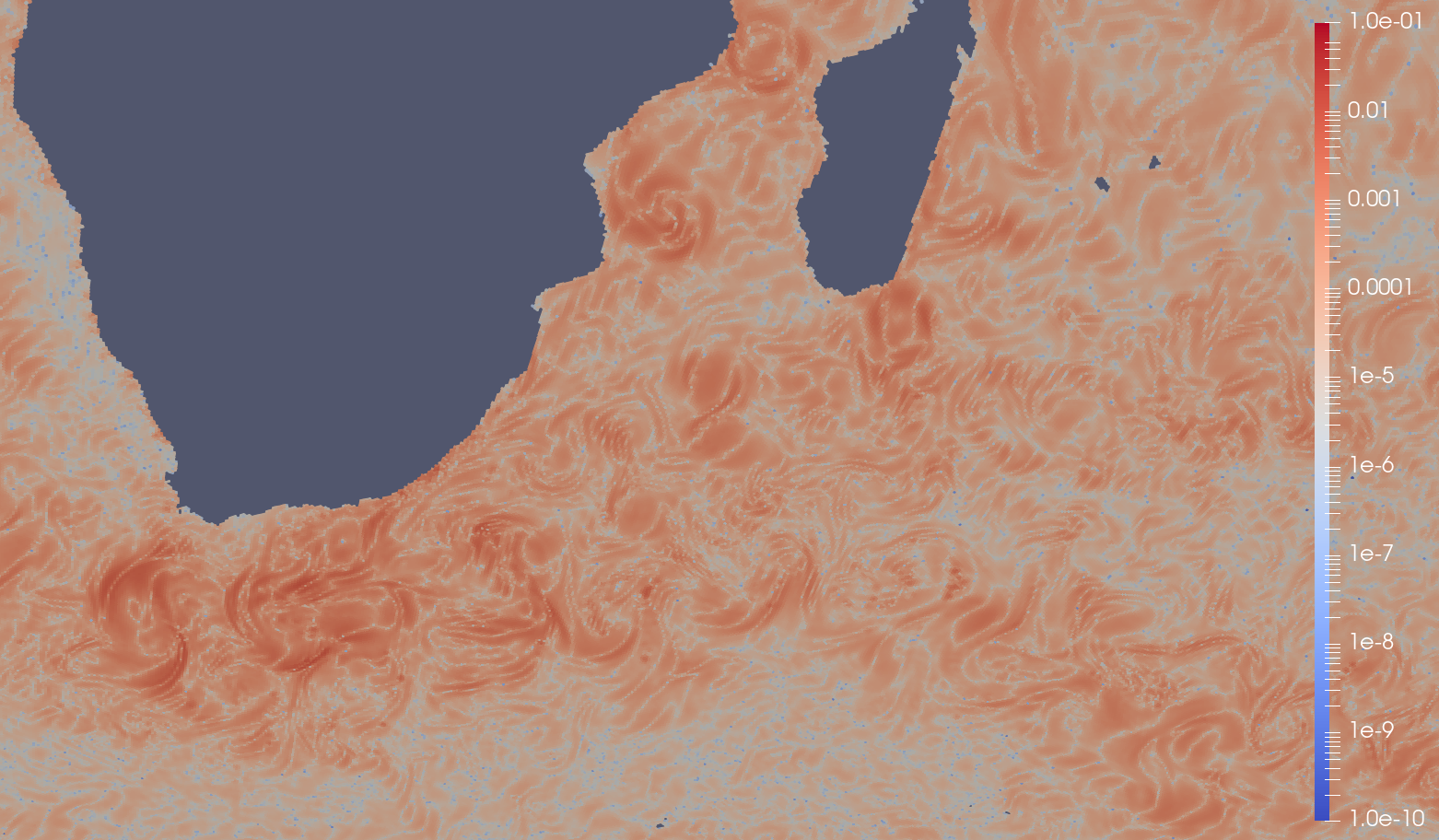}
    \caption{Figures showing magnitude of the difference of velocities for reference solution and simulated solution ($\abs{\bu_ref - \bu}$) measured on the ocean surface along southern coast of Africa on days 1 (left), 10 (middle), and 30 (right) for high resolution ``optimal'' trial.}
    \label{fig:abs error}
\end{figure}

\subsection{Minimum Length Scale}

One parameter of note was the length scale of the observational data, $\delta$. Here $\delta$ refers to the distance between observations measured in the number of cell lengths apart. The $\delta$ values were varied between $0$ and $3$ cell lengths, with a value of $0$ indicating that observations were taken from every cell in the domain. As one can see from \Cref{fig:h}, we obtained what looks like initial exponential decay, and then a roughly constant error for each value of $\delta = 0,1,2$. In the case of $\delta = 3$ we see slower exponential decay that is still present at day $30$ of our simulations; however, we expect similarly to our other $\delta$ values, that the error will level off as time increases.

\begin{figure}
    \centering
    \includegraphics[width=.8\textwidth]{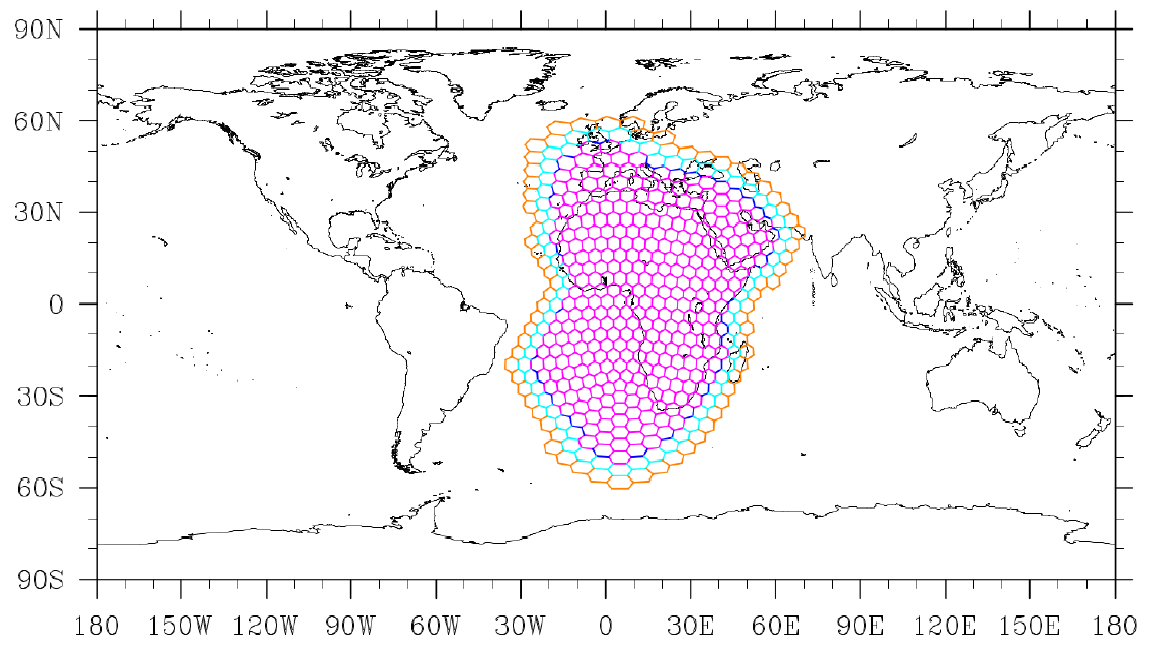}
    \caption{Representation of the edges owned by a single node. Pink edges are owned by the block itself, with the other colors representing edges in the halo layer that is shared with neighboring subdomain. The coloring (dark blue, light blue, and orange) correspond to the halo layer (first, second, and third), respectively.  (Image credit: MPAS Developer's guide \cite{mpas2013developers}.)}
    \label{fig:HALO}
\end{figure}
We note that the maximum value of $\delta = 3$ was chosen due to the limitations of the domain decomposition utilized for efficient parallel computation. Specifically, the domain is decomposed into smaller subdomains with an outer boundary ring, known as the \textit{halo}, that is $3$ cells in width, see \Cref{fig:HALO}. This halo is used to limit the communication between nodes by allowing each node to accurately access the information in the parts of the adjacent subdomains internally. That is, the interior of the domain is evolved forward in time using information from within the subdomain and from neighboring subdomains that is contained in the halo. Extending the value of $\delta$ past $3$ would require performing a halo update to obtain updated information from the boundary of the halo from each neighboring subdomain. Without performing such an update we cannot guarantee that the information contained in each node about the neighboring subdomains is accurate, or that our results will not vary between different runs of the same simulation. 

\subsection{Assimilation of tracer quantities}
Now we look into the effect of assimilating the tracer quantities of temperature and salinity. Due to the coupling between \Cref{momentum v} and \Cref{tracers v},  occurring through the pressure term, it is expected that assimilation of the tracer values may be required to obtain convergence to arbitrary precision. In this section we detail the results from simulations were we assimilated the tracer data in addition to velocity measurements.

We note here that in addition to the feedback control term on \Cref{momentum v} there is an additional feedback control term on \Cref{tracers v} to consider. As we stated in the preceding sections we note again that the parameter $\mu$ should have separate components $\mu_\bu$ and $\mu_{\phi}$ as the parameters should be adjusted based on the specific term they are applied to. In testing, we noted that the optimal choice of $\mu$ for the momentum equation was also the optimal choice for the tracer equations. Thus we opted to denote all of the feedback control parameters for each equation by a single $\mu$ parameter.

\begin{center}
    \begin{figure}
        \centering
        \includegraphics[width=12cm,height=8cm]{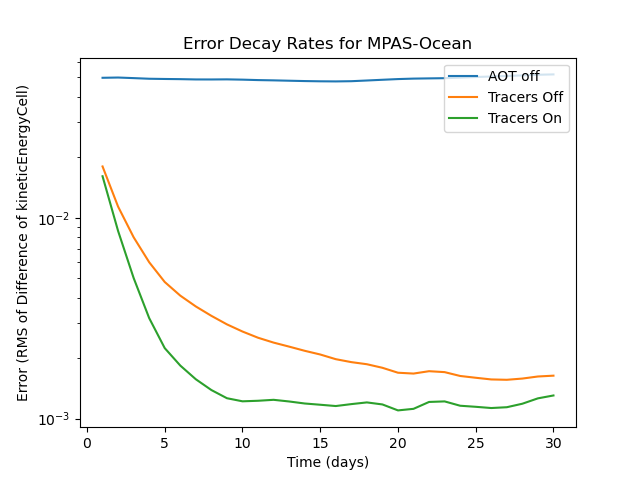}
        \caption{Error (in kinetic energy measurements) corresponding to the inclusion or exclusion of tracer based values.}
        \label{fig:tracers vel}
    \end{figure}
\end{center}

\begin{center}
     \begin{figure}
        \centering
\includegraphics[width=12cm,height=8cm]{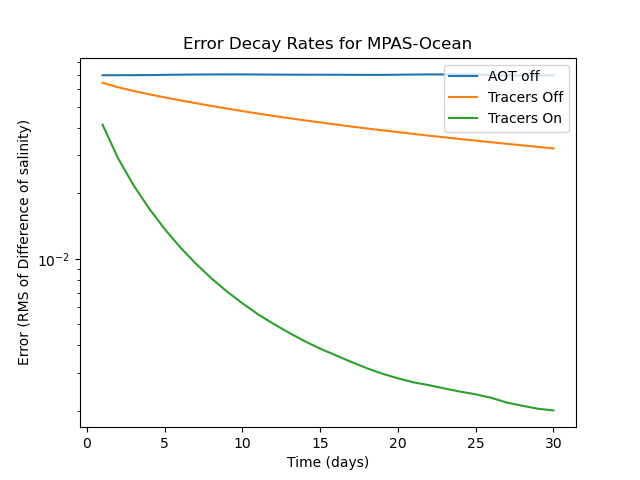}         \caption{Error (in salinity measurements) corresponding to the inclusion or exclusion of tracer based values.}
         \label{fig:tracers sal}
     \end{figure}
 \end{center}

\begin{center}
     \begin{figure}
        \centering
\includegraphics[width=12cm,height=8cm]{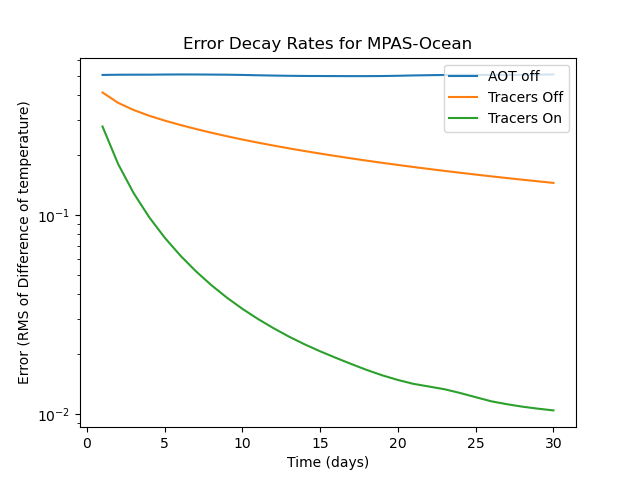} 
         \caption{Error (in temperature measurements) corresponding to the inclusion or exclusion of tracer based values.}
         \label{fig:tracers temp}
     \end{figure}
 \end{center}

As we see in \Cref{fig:tracers sal} the inclusion of tracer component data into the AOT algorithm does have an affect on convergence. However, the effect is still less than expected. 
In our simulations we were only able to obtain convergence to $1e-4$ with the addition of tracers. 
While the error does not decay to machine precision, we do observe several important features of the algorithm through the error decay. 
Namely, we see that the error decays exponentially fast in time in \emph{all} of the assimilated quantities: see, e.g., \Cref{fig:tracers vel,,fig:tracers temp,,fig:tracers sal}.

As the error only decays to the level of $1e-4$, there is the question of where this error is concentrated. It may be expected that when applying data assimilation on a realistic ocean model that error should accumulate on the coastlines and other domain boundaries. To clarify this, we see in \Cref{fig:abs error} that the absolute difference between observational data and the assimilated solution appears to be widespread and not restricted to coastlines. Moreover we see that the errors on day 10 appear to be qualitatively similar to those on day 30, so the error appears to converge to $1e-4$ throughout the domain with little decay past this level occurring.

\section{Conclusion}\label{conclusion}
In the present work, we have tested implementations of several versions of the Azouani-Olson-Titi (AOT) algorithm in the context of the MPAS-Ocean implementation of the primitive equations of the ocean.  We have observed that the AOT algorithm decreases the error exponentially fast before reaching a certain error level, dependent on the terms involved and whether the AOT feedback control term was handle implicitly or explicitly.
The range of error levels we observed in \Cref{table: schemes} covered roughly from $1e-4$ all the way to machine precision.

The reasons for this wide range are presently unknown.  From the table in \Cref{fig:valid dt}, we see that while including more parameters seems to increase the error, the Coriolis force, vertical mixing, and pressure gradient terms appear  to be the most problematic, as these terms each impart an order $1e-4$ error.  In the case of vertical mixing, we speculate that this may be due to the nature of the primitive equations themselves.  In particular, in the primitive equations, there is no evolution equation for the vertical velocity.  It has been observed by the authors in unpublished simulations of the other equations that when there is a prognostic variable (e.g., one without an evolution equation), or when the system is partially hyperbolic (say, a parabolic or diffusive system coupled to a hyperbolic or conservative system), the error can still decay, but tends to be bounded away from machine precision.

Moreover, in \cite{Pei_2019}, a lower bound on $\mu$ was given
in terms of a high power of the Grash\"of number. For ocean simulations the Grash\"of number is measured in millions, and so our range of $\mu$ values, constrained by a CFL condition, is not large enough in this setting to guarantee convergence according to the analytical theory. While numerical tests of the algorithm often indicate that such constraints are likely overly conservative by many orders of magnitude (see, e.g., \cite{Biswas_Bradshaw_Jolly_2021,Franz_Larios_Victor_2022,Gesho_Olson_Titi_2015,Larios_Pei_Victor_2023_second_best}), we cannot rule out the possibility that our parameter ranges are insufficient for machine-precision convergence to be realized.

While we did not typically observe convergence to machine precision, we note that our results appear to be consistent with the error levels observed in other studies applying the AOT algorithm in similar settings. In particular we note that in \cite{Carlson_VanRoekel_Petersen_Godinez_Larios_2021} 
the error levels were not seen to converge to arbitrary precision. We do note that in this previous study, utilization of the full solution as observational data was not considered, thus convergence to machine precision might not be expected. Regardless, we note that the convergence levels appear comparable and we are unaware of any studies at this time that have obtained convergence to machine precision for a high resolution ocean model using the AOT algorithm or a modification of it.

\section*{Acknowledgments}
 \noindent
This material is based upon work supported by the U.S. Department of Energy, Office of Science, Office of Workforce Development for Teachers and Scientists, Office of Science Graduate Student Research (SCGSR) program. The SCGSR program is administered by the Oak Ridge Institute for Science and Education (ORISE) for the DOE. ORISE is managed by ORAU under contract number DE-SC0014664. All opinions expressed in this paper are the authors' and do not necessarily reflect the policies and views of DOE, ORAU, or ORISE.

This research used resources of the National Energy Research Scientific Computing Center (NERSC), a U.S. DOE Office of Science User Facility located at Lawrence Berkeley National Laboratory, operated under Contract No. DE-AC02-05CH11231, as well as resources provided by the Los Alamos National Laboratory Institutional Computing Program, which is supported by the U.S. DOE National Nuclear Security Administration under Contract No. 89233218CNA000001.

The research of A.L. was supported in part by NSF grants DMS-2206762 and CMMI-1953346, and USGS  grant G23AS00157 number GRANT13798170.  
  The research of C.V. was supported in part by the NSF GRFP grant DMS-1610400.
M.R.P. was supported as part of the Energy Exascale Earth System Model (E3SM) project, funded by the U.S. Department of Energy, Office of Science, Office of Biological and Environmental Research.

\begin{scriptsize}
 \bibliographystyle{abbrv} 
 \bibliography{LariosBiblio.bib}

\begin{thebibliography}{10}

\bibitem{Albanez_Nussenzveig_Lopes_Titi_2016}
D.~A. Albanez, H.~J. Nussenzveig~Lopes, and E.~S. Titi.
\newblock Continuous data assimilation for the three-dimensional
  {N}avier--{S}tokes-$\alpha$ model.
\newblock {\em Asymptotic Anal.}, 97(1-2):139--164, 2016.

\bibitem{Anthes_1974_JAS}
R.~A. Anthes.
\newblock Data assimilation and initialization of hurricane prediction models.
\newblock {\em J. Atmos. Sci.}, 31(3):702--719, 1974.

\bibitem{Asch_Bocquet_Nodet_2016_DA_book}
M.~Asch, M.~Bocquet, and M.~Nodet.
\newblock {\em Data {A}ssimilation: {M}ethods, {A}lgorithms, and
  {A}pplications}.
\newblock SIAM, 2016.

\bibitem{Azouani_Olson_Titi_2014}
A.~Azouani, E.~Olson, and E.~S. Titi.
\newblock Continuous data assimilation using general interpolant observables.
\newblock {\em J. Nonlinear Sci.}, 24(2):277--304, 2014.

\bibitem{Azouani_Titi_2014}
A.~Azouani and E.~S. Titi.
\newblock Feedback control of nonlinear dissipative systems by finite
  determining parameters---a reaction-diffusion paradigm.
\newblock {\em Evol. Equ. Control Theory}, 3(4):579--594, 2014.

\bibitem{Bessaih_Olson_Titi_2015}
H.~Bessaih, E.~Olson, and E.~S. Titi.
\newblock Continuous data assimilation with stochastically noisy data.
\newblock {\em Nonlinearity}, 28(3):729--753, 2015.

\bibitem{Biswas_Bradshaw_Jolly_2022}
A.~Biswas, Z.~Bradshaw, and M.~Jolly.
\newblock Convergence of a mobile data assimilation scheme for the 2{D}
  {N}avier--{S}tokes equations.
\newblock {\em ar{X}iv:2210.11282}, 2022.

\bibitem{Biswas_Bradshaw_Jolly_2021}
A.~Biswas, Z.~Bradshaw, and M.~S. Jolly.
\newblock Data assimilation for the {N}avier--{S}tokes equations using local
  observables.
\newblock {\em SIAM J. Appl. Dyn. Syst.}, 20(4):2174--2203, 2021.

\bibitem{Carlson_Hudson_Larios_2020}
E.~Carlson, J.~Hudson, and A.~Larios.
\newblock Parameter recovery for the 2 dimensional {N}avier--{S}tokes equations
  via continuous data assimilation.
\newblock {\em SIAM J. Sci. Comput.}, 42(1):A250--A270, 2020.

\bibitem{Carlson_Hudson_Larios_Martinez_Ng_Whitehead_2021}
E.~Carlson, J.~Hudson, A.~Larios, V.~R. Martinez, E.~Ng, and J.~Whitehead.
\newblock Dynamically learning the parameters of a chaotic system using partial
  observations.
\newblock {\em Discrete Contin Dyn Syst Ser A}, 42(8):3809--3839, 2022.

\bibitem{Carlson_Larios_Titi_2023_nlDA}
E.~Carlson, A.~Larios, and E.~S. Titi.
\newblock Super-exponential convergence rate of a nonlinear continuous data
  assimilation algorithm: The 2{D} {N}avier--{S}tokes equations paradigm.
\newblock 2023.
\newblock (submitted) ar{X}iv:2304.01128.

\bibitem{Carlson_VanRoekel_Petersen_Godinez_Larios_2021}
E.~Carlson, L.~Van~Roekel, M.~Petersen, H.~C. Godinez, and A.~Larios.
\newblock {CDA} algorithm implemented in {MPAS-O} to improve eddy effects in a
  mesoscale simulation.
\newblock 2023.
\newblock (submitted).

\bibitem{Celik_Olson_2022}
E.~Celik and E.~Olson.
\newblock Data assimilation using time-delay nudging in the presence of
  {G}aussian noise.
\newblock {\em ar{X}iv e-prints}, 2022.

\bibitem{Celik_Olson_Titi_2019}
E.~Celik, E.~Olson, and E.~S. Titi.
\newblock Spectral filtering of interpolant observables for a discrete-in-time
  downscaling data assimilation algorithm.
\newblock {\em SIAM J. Appl. Dyn. Syst.}, 18(2):1118--1142, 2019.

\bibitem{Desamsetti_Dasari_Langodan_Knio_Hoteit_Titi_2019_WRF}
S.~Desamsetti, H.~Dasari, S.~Langodan, O.~Knio, I.~Hoteit, and E.~S. Titi.
\newblock Efficient dynamical downscaling of general circulation models using
  continuous data assimilation.
\newblock {\em Quarterly Journal of the Royal Meteorological Society}, 2019.

\bibitem{Diegel_Rebholz_2021}
A.~E. Diegel and L.~G. Rebholz.
\newblock Continuous data assimilation and long-time accuracy in a {$\rm C^0$}
  interior penalty method for the {C}ahn-{H}illiard equation.
\newblock {\em Appl. Math. Comput.}, 424:Paper No. 127042, 22, 2022.

\bibitem{Du_Shiue_2021}
Y.~J. Du and M.-C. Shiue.
\newblock Analysis and computation of continuous data assimilation algorithms
  for {L}orenz 63 system based on nonlinear nudging techniques.
\newblock {\em J. Computat. and Appl. Math.}, 386:113246, 2021.

\bibitem{Engwirda_2014}
D.~Engwirda.
\newblock {\em Locally-optimal Delaunay-refinement and optimisation-based mesh
  generation}.
\newblock PhD thesis, School of Mathematics and Statistics, The University of
  Sydney, 2014.

\bibitem{Engwirda_2015}
D.~Engwirda.
\newblock Voronoi-based point-placement for three-dimensional
  {D}elaunay-refinement.
\newblock {\em Procedia Engineering}, 124:330--342, 2015.

\bibitem{Engwirda_2016}
D.~Engwirda.
\newblock Conforming restricted delaunay mesh generation for piecewise smooth
  complexes.
\newblock {\em Procedia Engineering}, 163:84--96, 2016.

\bibitem{Engwirda_2018}
D.~Engwirda.
\newblock Generalised primal-dual grids for unstructured co-volume schemes.
\newblock {\em J. Comp. Phys.}, 375:155--176, 2018.

\bibitem{Engwirda_Ivers_2016}
D.~Engwirda and D.~Ivers.
\newblock Off-centre steiner points for delaunay-refinement on curved surfaces.
\newblock {\em Computer-Aided Design}, 72:157--171, 2016.

\bibitem{Farhat_GlattHoltz_Martinez_McQuarrie_Whitehead_2019}
A.~Farhat, N.~E. Glatt-Holtz, V.~R. Martinez, S.~A. McQuarrie, and J.~P.
  Whitehead.
\newblock Data {A}ssimilation in {L}arge {P}randtl {R}ayleigh--{B}\'{e}nard
  {C}onvection from {T}hermal {M}easurements.
\newblock {\em SIAM J. Appl. Dyn. Syst.}, 19(1):510--540, 2020.

\bibitem{Farhat_Jolly_Titi_2015}
A.~Farhat, M.~S. Jolly, and E.~S. Titi.
\newblock Continuous data assimilation for the 2{D} {B}\'enard convection
  through velocity measurements alone.
\newblock {\em Phys. D}, 303:59--66, 2015.

\bibitem{Farhat_Larios_Martinez_Whitehead_2023_force}
A.~Farhat, A.~Larios, V.~R. Martinez, and J.~P. Whitehead.
\newblock Identifying the body force from partial observations of a 2{D}
  incompressible velocity field.
\newblock {\em (submitted) ar{X}iv:2302.04701}, 2023.

\bibitem{Farhat_Lunasin_Titi_2016abridged}
A.~Farhat, E.~Lunasin, and E.~S. Titi.
\newblock Abridged continuous data assimilation for the 2{D} {N}avier--{S}tokes
  equations utilizing measurements of only one component of the velocity field.
\newblock {\em J. Math. Fluid Mech.}, 18(1):1--23, 2016.

\bibitem{Farhat_Lunasin_Titi_2016benard}
A.~Farhat, E.~Lunasin, and E.~S. Titi.
\newblock Data assimilation algorithm for 3{D} {B\'e}nard convection in porous
  media employing only temperature measurements.
\newblock {\em J. Math. Anal. Appl.}, 438(1):492--506, 2016.

\bibitem{Farhat_Lunasin_Titi_2016_Charney}
A.~Farhat, E.~Lunasin, and E.~S. Titi.
\newblock On the {C}harney conjecture of data assimilation employing
  temperature measurements alone: the paradigm of 3{D} planetary geostrophic
  model.
\newblock {\em Mathematics of Climate and Weather Forecasting}, 2(1), 2016.

\bibitem{Farhat_Lunasin_Titi_2017_Horizontal}
A.~Farhat, E.~Lunasin, and E.~S. Titi.
\newblock Continuous data assimilation for a {2D} {B}\'enard convection system
  through horizontal velocity measurements alone.
\newblock {\em J. Nonlinear Sci.}, pages 1--23, 2017.

\bibitem{Foias_Mondaini_Titi_2016}
C.~Foias, C.~F. Mondaini, and E.~S. Titi.
\newblock A discrete data assimilation scheme for the solutions of the
  two-dimensional {N}avier--{S}tokes equations and their statistics.
\newblock {\em SIAM J. Appl. Dyn. Syst.}, 15(4):2109--2142, 2016.

\bibitem{Franz_Larios_Victor_2022}
T.~Franz, A.~Larios, and C.~Victor.
\newblock The bleeps, the sweeps, and the creeps: {C}onvergence rates for
  dynamic observer patterns via data assimilation for the 2{D}
  {N}avier-{S}tokes equations.
\newblock {\em Comput. Methods Appl. Mech. Engrg.}, 392:Paper No. 114673, 19,
  2022.

\bibitem{Gardner_Larios_Rebholz_Vargun_Zerfas_2020_VVDA}
M.~Gardner, A.~Larios, L.~G. Rebholz, D.~Vargun, and C.~Zerfas.
\newblock Continuous data assimilation applied to a velocity-vorticity
  formulation of the 2{D} {N}avier--{S}tokes equations.
\newblock {\em Electron. Res. Arch.}, 29(3):2223--2247, 2021.

\bibitem{Gesho_Olson_Titi_2015}
M.~Gesho, E.~Olson, and E.~S. Titi.
\newblock A computational study of a data assimilation algorithm for the
  two-dimensional {N}avier--{S}tokes equations.
\newblock {\em Commun. Comput. Phys.}, 19(4):1094--1110, 2016.

\bibitem{mpas-3}
J.-C. Golaz, P.~M. Caldwell, L.~P.~V. Roekel, M.~R. Petersen, Q.~Tang, J.~D.
  Wolfe, G.~Abeshu, V.~Anantharaj, X.~S. Asay-Davis, D.~C. Bader, S.~A.
  Baldwin, G.~Bisht, P.~A. Bogenschutz, M.~Branstetter, M.~A. Brunke, S.~R.
  Brus, S.~M. Burrows, P.~J. Cameron-Smith, A.~S. Donahue, M.~Deakin, R.~C.
  Easter, K.~J. Evans, Y.~Feng, M.~Flanner, J.~G. Foucar, J.~G. Fyke, B.~M.
  Griffin, C.~Hannay, B.~E. Harrop, M.~J. Hoffman, E.~C. Hunke, R.~L. Jacob,
  D.~W. Jacobsen, N.~Jeffery, P.~W. Jones, N.~D. Keen, S.~A. Klein, V.~E.
  Larson, L.~R. Leung, H.-Y. Li, W.~Lin, W.~H. Lipscomb, P.-L. Ma, S.~Mahajan,
  M.~E. Maltrud, A.~Mametjanov, J.~L. McClean, R.~B. McCoy, R.~B. Neale, S.~F.
  Price, Y.~Qian, P.~J. Rasch, J.~E. J.~R. Eyre, W.~J. Riley, T.~D. Ringler,
  A.~F. Roberts, E.~L. Roesler, A.~G. Salinger, Z.~Shaheen, X.~Shi, B.~Singh,
  J.~Tang, M.~A. Taylor, P.~E. Thornton, A.~K. Turner, M.~Veneziani, H.~Wan,
  H.~Wang, S.~Wang, D.~N. Williams, P.~J. Wolfram, P.~H. Worley, S.~Xie,
  Y.~Yang, J.-H. Yoon, M.~D. Zelinka, C.~S. Zender, X.~Zeng, C.~Zhang,
  K.~Zhang, Y.~Zhang, X.~Zheng, T.~Zhou, and Q.~Zhu.
\newblock The {DOE E3SM} coupled model version 1: {O}verview and evaluation at
  standard resolution.
\newblock {\em J. Adv. Model. Earth Sy.}, 11(7):2089--2129, 2019.

\bibitem{Hammoud_Titi_Hoteit_Knio_2022}
M.~A. E.~R. Hammoud, E.~S. Titi, I.~Hoteit, and O.~Knio.
\newblock Cdanet: A physics-informed deep neural network for downscaling fluid
  flows.
\newblock {\em Journal of Advances in Modeling Earth Systems},
  14(12):e2022MS003051, 2022.

\bibitem{Hayden_Olson_Titi_2011}
K.~Hayden, E.~Olson, and E.~S. Titi.
\newblock Discrete data assimilation in the {L}orenz and 2{D}
  {N}avier--{S}tokes equations.
\newblock {\em Phys. D}, 240(18):1416--1425, 2011.

\bibitem{Hoke_Anthes_1976_MWR}
J.~E. Hoke and R.~A. Anthes.
\newblock The initialization of numerical models by a dynamic-initialization
  technique.
\newblock {\em Monthly Weather Review}, 104(12):1551--1556, 1976.

\bibitem{Jolly_Martinez_Olson_Titi_2018_blurred_SQG}
M.~S. Jolly, V.~R. Martinez, E.~J. Olson, and E.~S. Titi.
\newblock Continuous data assimilation with blurred-in-time measurements of the
  surface quasi-geostrophic equation.
\newblock {\em Chin. Ann. Math. Ser. B}, 40(5):721--764, 2019.

\bibitem{Jolly_Martinez_Titi_2017}
M.~S. Jolly, V.~R. Martinez, and E.~S. Titi.
\newblock A data assimilation algorithm for the subcritical surface
  quasi-geostrophic equation.
\newblock {\em Adv. Nonlinear Stud.}, 17(1):167--192, 2017.

\bibitem{Lakshmivarahan_Lewis_2013}
S.~Lakshmivarahan and J.~M. Lewis.
\newblock Nudging methods: {A} critical overview.
\newblock In {\em Data Assimilation for Atmospheric, Oceanic and Hydrologic
  Applications (Vol. II)}, pages 27--57. Springer, 2013.

\bibitem{Larios_Pei_2018_NSV_DA}
A.~Larios and Y.~Pei.
\newblock Approximate continuous data assimilation of the 2{D}
  {N}avier--{S}tokes equations via the {V}oigt-regularization with observable
  data.
\newblock {\em Evol. Equ. Control Theory}, 9(3):733--751, 2020.

\bibitem{Larios_Pei_2017_KSE_DA_NL}
A.~Larios and Y.~Pei.
\newblock Nonlinear continuous data assimilation.
\newblock {\em Evol. Equ. Control Theory}, 2024.
\newblock (accepted for publication).

\bibitem{Larios_Rebholz_Zerfas_2018}
A.~Larios, L.~G. Rebholz, and C.~Zerfas.
\newblock Global in time stability and accuracy of {IMEX}-{FEM} data
  assimilation schemes for {N}avier--{S}tokes equations.
\newblock {\em Comput. Methods Appl. Mech. Engrg.}, 345:1077--1093, 2019.

\bibitem{Larios_Victor_2019}
A.~Larios and C.~Victor.
\newblock Continuous data assimilation with a moving cluster of data points for
  a reaction diffusion equation: a computational study.
\newblock {\em Commun. Comput. Phys.}, 29(4):1273--1298, 2021.

\bibitem{Larios_Pei_Victor_2023_second_best}
A.~Larios and C.~Victor.
\newblock The second-best way to do sparse-in-time continuous data
  assimilation: {I}mproving convergence rates for the 2{D} and 3{D}
  {N}avier--{S}tokes equations.
\newblock {\em (submitted) ar{X}iv:2303.03495}, 2023.

\bibitem{Law_Stuart_Zygalakis_2015_book}
K.~Law, A.~Stuart, and K.~Zygalakis.
\newblock {\em A {M}athematical {I}ntroduction to {D}ata {A}ssimilation},
  volume~62 of {\em Texts in Applied Mathematics}.
\newblock Springer, Cham, 2015.

\bibitem{Lunasin_Titi_2015}
E.~Lunasin and E.~S. Titi.
\newblock Finite determining parameters feedback control for distributed
  nonlinear dissipative systems--a computational study.
\newblock {\em Evol. Equ. Control Theory}, 6(4):535--557, 2017.

\bibitem{Martinez_2022}
V.~R. Martinez.
\newblock Convergence analysis of a viscosity parameter recovery algorithm for
  the 2{D} {N}avier-{S}tokes equations.
\newblock {\em Nonlinearity}, 35(5):2241--2287, 2022.

\bibitem{mpas2013developers}
{MPAS Development Team}.
\newblock Mpas developers guide.
\newblock Online, Nov. 2013.

\bibitem{Olson_Titi_2003}
E.~Olson and E.~S. Titi.
\newblock Determining modes for continuous data assimilation in 2{D}
  turbulence.
\newblock {\em J. Statist. Phys.}, 113(5-6):799--840, 2003.
\newblock Progress in statistical hydrodynamics (Santa Fe, NM, 2002).

\bibitem{Olson_Titi_2008_TCFD}
E.~Olson and E.~S. Titi.
\newblock Determining modes and {G}rashof number in 2{D} turbulence: a
  numerical case study.
\newblock {\em Theor. Comp. Fluid Dyn.}, 22(5):327--339, 2008.

\bibitem{Pachev_Whitehead_McQuarrie_2021concurrent}
B.~Pachev, J.~P. Whitehead, and S.~A. McQuarrie.
\newblock Concurrent multiparameter learning demonstrated on the
  {K}uramoto--{S}ivashinsky equation.
\newblock {\em SIAM J. Sci. Comput.}, 44(5):A2974--A2990, 2022.

\bibitem{Pei_2019}
Y.~Pei.
\newblock Continuous data assimilation for the 3{D} primitive equations of the
  ocean.
\newblock {\em Commun. Pure Appl. Anal.}, 18(2):643--661, 2019.

\bibitem{mpasoceanuser}
M.~R. Petersen, X.~S. Asay-Davis, D.~W. Jacobsen, M.~E. Maltrud, T.~D. Ringler,
  L.~Van~Roekel, C.~Veneziani, and P.~J. Wolfram~Jr.
\newblock Mpas-ocean model user's guide version 6.0.
\newblock Technical report, Los Alamos National Lab.(LANL), Los Alamos, NM
  (United States), 2018.

\bibitem{mpas-1}
T.~Ringler, M.~Petersen, R.~L. Higdon, D.~Jacobsen, P.~W. Jones, and
  M.~Maltrud.
\newblock A multi-resolution approach to global ocean modeling.
\newblock {\em Ocean Modelling}, 69:211--232, 2013.

\bibitem{yu2019twin}
L.~Yu, K.~Fennel, B.~Wang, A.~Laurent, K.~R. Thompson, and L.~K. Shay.
\newblock Evaluation of nonidentical versus identical twin approaches for
  observation impact assessments: an ensemble-{K}alman-filter-based ocean
  assimilation application for the {G}ulf of {M}exico.
\newblock {\em Ocean Science}, 15(6):1801--1814, 2019.

\end{thebibliography}
 \end{scriptsize}





\end{document}